\documentclass[11pt]{article}
\usepackage{graphics,graphicx}
\usepackage{amssymb,amsmath}
\usepackage{enumitem}

\usepackage[paperwidth=216mm, paperheight=279mm, margin=25.4mm]{geometry}

\newtheorem{thm}{Theorem}

\newtheorem{problem}{Problem}
\newtheorem{conj}{Conjecture}
\usepackage{thmtools}
\newlist{lemmalist}{enumerate}{1}
\setlist[lemmalist]{label=(\roman{lemmalisti}),
	ref=\thelemma:$(\roman{lemmalisti})$,
	noitemsep}
\declaretheorem[
name=Lemma]{lemma}
\usepackage[capitalize]{cleveref}

\Crefname{lemmalist}{Lemma}{Lemmas}
\addtotheorempostheadhook[lemma]{\crefalias{lemmalisti}{lemmalist}}

\newtheorem{pro}{}[thm]

\newtheorem{claim}{Claim}[lemma]

\newtheorem{cla}{Claim}[thm]

\newtheorem{prop}{}

\def \no {\noindent}
 
 \def \sm {\setminus}
 \def \es {\emptyset}

\newcommand*\sq{\mathbin{\vcenter{\hbox{\rule{0.75ex}{0.75ex}}}}}

\begin{document}

\title{Coloring of ($P_5$, 4-wheel)-free graphs}
\author{ Arnab Char\thanks{Computer Science Unit, Indian Statistical
Institute, Chennai Centre, Chennai 600029, India.} \and T.~Karthick\thanks{Corresponding author, Computer Science Unit, Indian Statistical
Institute, Chennai Centre, Chennai 600029, India. Partially supported by DST-SERB, Government of India, under MATRICS scheme (MTR/2018/000288).}}

\date{\today}

\maketitle

\begin{abstract}
For a graph $G$,  $\chi(G)$ $(\omega(G))$  denote its chromatic (clique) number. A $P_5$ is the chordless path on five vertices, and a \emph{$4$-wheel}  is  the graph
consisting of a chordless cycle on four vertices $C_4$ plus an additional vertex adjacent to all the vertices of the $C_4$. In this paper, we show that  every ($P_5$,\,$4$-wheel)-free graph $G$  satisfies $\chi(G)\leq \frac{3}{2}\omega(G)$.  Moreover,  this bound is almost tight. That is, there is a class of  ($P_5$,\,$4$-wheel)-free graphs $\cal L$ such that every graph $H\in \cal L$ satisfies $\chi(H)\geq\frac{10}{7}\omega(H)$. This generalizes/improves several previously known results in the literature.
\end{abstract}

\section{Introduction}
All our graphs are simple and finite.
Given a graph $G$, as usual, we write $\chi(G)$ to denote the chromatic number of $G$,  and $\omega(G)$ to denote the size of a maximum  clique in  $G$.
A graph $G$ is \emph{perfect}, if every induced subgraph $H$ of $G$ satisfies $\chi(H)=\omega(H)$.  As introduced by Gy\'arf\'as \cite{Gyarfas}, a hereditary class of graphs $\cal C$ is said to be \emph{$\chi$-bounded}, if there is a function $f:\mathbb{Z}^+\rightarrow \mathbb{Z}^+$ (called a \emph{$\chi$-binding function} for $\cal C$) such that every $G\in \cal C$ satisfies $\chi(G)\leq f(\omega(G))$. Obviously the class of perfect graphs is $\chi$-bounded with $f(x)=x$ as a $\chi$-binding function.  Recently there has been much research on $\chi$-bounded classes of graphs; see \cite{CHPV, CPST, survey, SS,  TP-prob}  for examples. We refer to \cite{ScottSey-Survey} for a comprehensive survey  on  $\chi$-bounded classes of graphs and their connections to other topics in graph theory.

Given a positive integer $k$, let $P_{k}$ denote the chordless (or induced) path on $k$ vertices, and for $k\ge 3$, $C_{k}$ denote the chordless (or induced) cycle on $k$ vertices. For $k\geq 4$, a \emph{$k$-wheel}  is  the graph
consisting of a cycle $C_k$  plus an additional vertex adjacent to all the vertices of the $C_k$. We say that a graph $G$
\emph{contains} a graph $H$, if $G$ has an induced subgraph isomorphic to $H$.
Given a family of graphs $\cal F$,  a graph $G$ is  $\cal F$-free if it does not contain any member of $\cal F$.

Gy\'arf\'as \cite{Gyarfas} showed that the class of $P_t$-free graphs is $\chi$-bounded with $\chi$-binding function $f(x)=(t-1)^{x-1}$. Gravier, Ho\'ang and Maffray \cite{GHM} improved this result, and
showed that for $t\geq 4$ and $\omega(G)\geq 2$, every $P_t$-free graph $G$ satisfies $\chi(G)\leq (t-2)^{\omega(G)-1}$.
 It is well-known that every $P_4$-free graph is perfect.
Esperet,  Lemoine, Maffray, Morel \cite{ELMM} showed that  every $P_5$-free graph $G$ with $\omega(G)\geq 3$ satisfies $\chi(G)\leq 5\cdot3^{\omega(G)-3}$. Recently, Trotignon and Pham \cite{TP-prob} posed the following question  (see also \cite{ScottSey-Survey}):
 \begin{problem}  Is it true that,  the class of $P_5$-free graphs is polynomially $\chi$-bounded?
 \end{problem}
More generally, Esperet \cite{Esperet} conjectured the following:
\begin{conj}  If $\cal G$ is  a $\chi$-bounded class of graphs, then  $\cal G$ is polynomially $\chi$-bounded.
 \end{conj}

 Problem~1 is open even for a subclass of the class of $P_5$-free graphs, namely the class of ($P_5, C_5$)-free graphs. Chudnovsky and Sivaraman \cite{Chud-Siv} showed that every ($P_5,\,C_5$)-free graph $G$ satisfies $\chi(G)\leq 2^{\omega(G)-1}$. Moreover, Fouquet, Giakoumakis, Maire and Thuillier \cite{Fouquet} showed that there does not exist a  linear $\chi$-binding function even for the class of ($P_5$,\,$P_5^c$)-free graphs. It is interesting to note that the existence of a polynomial $\chi$-binding function for the class of $P_5$-free graphs implies the Erd\"os-Hajnal conjecture for the class of $P_5$-free graphs; see \cite{ScottSey-Survey}. Recently Scott, Seymour and Spirkl \cite{SSS} showed that if $G$ is $P_5$-free and $\omega(G)\geq  3$,  then $\chi(G) \leq \omega(G)^{\log_2(\omega(G))}$.

In this paper, we are interested in  (smallest) $\chi$-binding functions for the class of ($P_5$,\,$H$)-free graphs, where $H$ is a   connected graph on at most five vertices. Only very few graph classes have been investigated in this direction. It is well-known that every $P_4$-free graph is perfect.  Esperet et al. \cite{ELMM} showed that every  ($P_5$,\,$K_4$)-free graph $G$ satisfies $\chi(G)\leq 5$. It follows from a result of Randerath \cite{Rand-thesis} that  every  ($P_5$,\,paw)-free graph $G$ satisfies $\chi(G)\leq \omega(G)+1$.  Recently, the second author with Chudnovsky, Maceli  and Maffray \cite{CKMM-gem} showed that every ($P_5$,\,gem)-free graph $G$ satisfies $\chi(G)\leq \lceil\frac{5\omega(G)}{4}\rceil$, and with Huang \cite{KH-Paragliger}, he showed that every ($P_5$,\,paraglider)-free graph $G$ satisfies $\chi(G)\leq \frac{3}{2}\omega(G)$. Chudnovsky and Sivaraman \cite{Chud-Siv} showed that if a graph $G$ is ($P_5$, bull)-free, then $\chi(G)\leq \binom{\omega(G)+1}{2}$.    Shiermeyer \cite{Sch} recently studied the $\chi$-binding functions for  ($P_5$,\,butterfly)-free graphs and ($P_5$,\,hammer)-free graphs. In this paper, we study the class of  ($P_5$,\,$4$-wheel)-free graphs, and prove the following.

 \begin{thm}\label{thm:32-bnd} If $G$ is a ($P_5$,\,$4$-wheel)-free graph, then $\chi(G)\le \frac{3}{2}\omega(G)$. Moreover, there is a class of ($P_5$,\,$4$-wheel)-free graphs $\cal L$ such that every graph $H \in \cal L$ satisfies $\chi(H)\geq \frac{10}{7}\omega(H)$.
 \end{thm}

 We  note that Theorem~\ref{thm:32-bnd} generalizes/improves the following known results.
\begin{itemize}\itemsep=0pt
\item For every ($2K_2$,\,$C_4$)-free graph $G$, $\chi(G)$ is either $\omega(G)$ or $\omega(G)+1$  \cite{Tuza, Gyarfas}.
\item Every ($2K_2$,\,$4$-wheel)-free graph satisfies $\chi(G)\leq \omega(G)+5$ \cite{KSM-DM}.
\item Every ($3K_1$,\,$4$-wheel)-free graph satisfies $\chi(G)\leq 2\omega(G)$ \cite{CKS2008}.
\item Every ($P_5$, diamond)-free graph satisfies $\chi(G)\leq \omega(G)+1$ \cite{survey}.
\item Every ($P_5$,\,$C_4$)-free graph $G$ satisfies $\chi(G)\leq \lceil\frac{5\omega(G)}{4}\rceil$ \cite{CKS2007}.
\item Every ($P_5$,\,$4$-wheel)-free graph $G$ satisfies $\chi(G)\leq 5\lceil\frac{5\omega(G)}{4}\rceil$ \cite{CKS2007}.
\end{itemize}

The proof of Theorem~\ref{thm:32-bnd} is derived from the structure theorem for the class of ($P_5$, $4$-wheel)-free graphs (Theorem~\ref{thm:struc}) given below. Before stating it, we present some definitions.
A \emph{clique} (\emph{stable set}) in $G$ is a set of pairwise adjacent (nonadjacent) vertices in $G$. A \emph{clique cutset} in   $G$ is a clique  $Q$ in $G$  such that $G- Q$ has more connected components than $G$. A graph is an \emph{atom} if it has no clique cutset.  A graph $G$ is a \emph{quasi-line graph} if for each $v\in V(G)$, the set of neighbors of $v$ can be expressed as the union of two cliques.

For any two sets $X, Y \subset V(G)$ of a graph $G$, we say that $X$ is complete (anticomplete) to $Y$
if every vertex in $X$ is adjacent (nonadjacent) to every vertex in $Y$. If $X$ is singleton, say $X=\{v\}$, then we simply  write $v$ is complete (anticomplete) to $Y$, instead of  $\{v\}$ is complete (anticomplete) to $Y$. Given two vertex-subsets $X$ and $Y$, we say that $X$ \emph{meets} $Y$ if $X\cap Y\neq \es$.

\begin{figure}[t]
\centering
 \includegraphics[height=6cm]{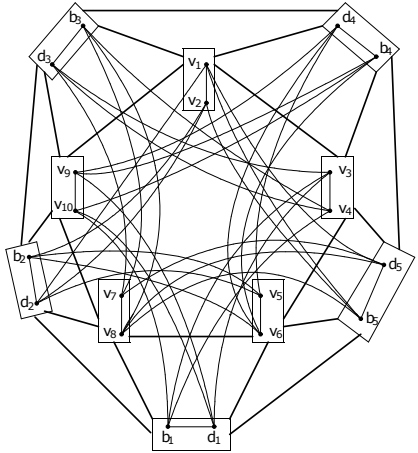}
\caption{Example of a ($P_5$,\,$4$-wheel)-free graph $G^*$ with chromatic number 10 and clique number~7 (see also \cite{CKS2008}). Here, a bold (or thick) line between two rectangles represents that every vertex inside a rectangle is adjacent to every vertex in the other. For example, the vertex $v_5$ is adjacent to both $b_1$ and $d_1$. Likewise, $v_5$ is adjacent to both $v_7$ and $v_8$.}\label{fig:10by7-example}
\end{figure}

We say that a graph $G$ is  \emph{nice} if there are three disjoint stable sets $S_1, S_2$ and $S_3$ whose union meets each maximum clique of $G$ at least twice, and the other maximal cliques of $G$ at least once, that is, $\omega(G-(S_1\cup S_2\cup S_3)) \le \omega(G)-2$.

Now we can state our  structure theorem.
\begin{thm}\label{thm:struc}
If $G$ is a connected ($P_5$,\,$4$-wheel)-free atom, then  $G$ is either a perfect graph, a nice graph,  or   a quasi-line graph.
\end{thm}

Theorem~\ref{thm:struc} is derived from Theorem~\ref{thm:struc-new}  given below. (The proof  of Theorem~\ref{thm:struc-new} is given in Section~\ref{sec-mainstruc}.)

\begin{thm}\label{thm:struc-new}
If $G$ is a connected ($P_5$,\,$4$-wheel)-free  atom, then one of the following holds:
\begin{enumerate}
\item[(1)] If $G$ has an induced $5$-wheel, then $G$ is a nice graph.
\item[(2)] If $G$  has an induced $C_5$ and has no induced $5$-wheel, then  $G$ is either a nice graph or   a quasi-line graph.
\item[(3)] If $G$ has an induced $C_7^c$ and has no induced $C_5$, then   $G$ is a nice graph.
\end{enumerate}
\end{thm}

\medskip
\no{\it Proof of Theorem~\ref{thm:struc}, assuming Theorem~\ref{thm:struc-new}}. Let $G$ be a connected ($P_5$,\,$4$-wheel)-free atom. We may assume that $G$ is not perfect. Then since $C_{2k+1}$ for $k\geq 3$ contains an induced $P_5$,  and since $C_{2k+1}^c$ for $k\geq 4$ contains an induced $4$-wheel,
by the Strong Perfect Graph Theorem \cite{spgt}, $G$ contains an induced $C_5 (\cong C_5^c)$ or $C_{7}^c$. So it satisfies the hypothesis of one of the items of
Theorem~\ref{thm:struc-new} and subsequently it satisfies the conclusion of this item.   This proves Theorem~\ref{thm:struc}. \hfill{$\Box$}

\medskip
 Next we give a proof of Theorem~\ref{thm:32-bnd} using Theorem~\ref{thm:struc}. A \emph{blowup} of a graph $H$ is any graph $G$ such that $V(G)$ can be partitioned
into $|V(H)|$  (not necessarily nonempty) sets $Q_v$, $v\in  V(H)$, such that  each $Q_v$ induces a $P_3$-free graph,   $Q_u$ is complete to $Q_v$ if $uv\in E(H)$, and $Q_u$ is anticomplete to $Q_v$ if $uv\notin E(H)$. A blowup is a \emph{clique-blowup} if each  $Q_v$ is a clique.

\medskip
\no{\it Proof of Theorem~\ref{thm:32-bnd}}. Let $G$ be a ($P_5$, $4$-wheel)-free graph. We prove the first assertion by induction  on $|V(G)|$.   We may assume that $G$ is connected and not perfect.
If $G$ has a clique cutset, say $Q$, let $V_1, V_2$ be a partition of $V(G)\sm Q$ such that $V_1, V_2\neq\es$, and $V_1$ is anticomplete to $V_2$. Then $\chi(G) = \max\{\chi(G[Q\cup V_1]), \chi(G[Q\cup V_2])\}\leq \max\{\frac{3}{2}\omega(G[Q\cup V_1]), \frac{3}{2}\omega(G[Q\cup V_2])\} \leq \frac{3}{2}\omega(G)$.
So we may assume that $G$ is an atom, and we apply Theorem~\ref{thm:struc}. If $G$ is a quasi-line graph, then by a result of Chudnovsky and Ovetsky \cite{Chud-Quasiline}, we have  $\chi(G)\le  \frac{3}{2}\omega(G)$. So suppose that $G$ is nice. Then $G$ has three stable sets
$S_1,S_2,$ and $S_3$ such that $\omega(G- (S_1\cup S_2\cup S_3))\le
\omega(G)-2$. Consider any $\chi(G- (S_1\cup S_2\cup S_3))$-coloring of
$G- (S_1\cup S_2\cup S_3)$, and use $S_1, S_2,$ and $S_3$ as three new
color classes to get a coloring of $G$.  Then we have $\chi(G)\le \chi(G-
(S_1\cup S_2\cup S_3))+3\le \frac{3}{2}(\omega(G)-2)
+3 =\frac{3}{2}\omega(G)$.
 This proves the first assertion of Theorem~\ref{thm:32-bnd}.

To prove the  second assertion of Theorem~\ref{thm:32-bnd}, consider the clique-blowup $H$ of the graph $G^*$ shown in Figure~\ref{fig:10by7-example} where $|Q_v|=k \ge 1$, for each $v\in V(G^*)$. Then it is shown in \cite{CKS2008} that $H$ is ($3K_1$,\,$4$-wheel)-free (and hence ($P_5$,\,$4$-wheel)-free), and  that $\omega(H)=7k$. Moreover, since $H$ has no stable set of size $3$, $\chi(H)\geq \frac{|V(H)|}{2} = \frac{20k}{2} = 10k$. This completes the proof of Theorem~\ref{thm:32-bnd}.  \hfill{$\Box$}

\smallskip

The remainder of the paper is organized as follows.
In Section~\ref{sec-prel}, we give some preliminaries which are used in this paper. In Section~\ref{sec-C5}, we present some useful structural properties of ($P_5$,\,$4$-wheel)-free atoms that has an induced $C_5$, and finally in Section~\ref{sec-mainstruc}, we  prove Theorem~\ref{thm:struc-new}.

\section{Preliminaries}\label{sec-prel}
We follow West \cite{west} for standard notation  and terminology. Let $G$ be a graph with vertex-set
$V(G)$ and edge-set $E(G)$. The \emph{complement graph} of $G$ is denoted by $G^c$.   If $X \subseteq V(G)$, then $G[X]$ and $G-X$ respectively denote the subgraph induced by $X$ and $V(G)\sm X$ in $G$.    Given $u,v \in V(G)$, we say that a vertex $u$ is a neighbor of $v$ if $u$ and $v$ are adjacent in $G$. The \emph{neighborhood} of a vertex $v$, denoted by $N_G(v)$, is the set of neighbors of $v$ in $G$ (and we drop the
subscript $G$ when there is no ambiguity). If $X \subseteq V(G)$, then $N(X)$ denote the set $\{x\in V(G)\sm X\mid x \mbox{ has a neighbor in } X\}$.

Given two vertex-disjoint graphs $G$ and $H$, the
\emph{union} $G\cup H$, is the graph with vertex-set $V(G)\cup V(H)$ and edge-set $E(G)\cup E(H)$.  The union of $k$ graphs each isomorphic to $G$ is denoted by $kG$; for
instance $2K_2$ denotes the graph that consists union of two disjoint  $K_2$'s.  An induced cycle $C_k$ with vertex-set $\{v_1,v_2,\ldots, v_k\}$ and edge-set $\{v_1v_2,v_2v_3,\ldots, v_{k-1}v_k, v_kv_1\}$ will be simply denoted by
$v_1$-$v_2$-$\cdots$-$v_k$-$v_1$. Likewise, an induced path $P_k$ with vertex-set $\{v_1,v_2,\ldots, v_k\}$ and edge-set $\{v_1v_2,v_2v_3,\ldots, v_{k-1}v_k\}$ will be simply denoted by
$v_1$-$v_2$-$\cdots$-$v_k$. For $k\geq 4$, a \emph{$k$-wheel} is the graph with   vertex-set $\{v_1, v_2,\ldots, v_k, c\}$ such that $v_1$-$v_2$-$\cdots$-$v_k$-$v_1$ is a $C_k$, and $c$ is complete to $\{v_1, \ldots,v_k\}$. A graph is said to be \emph{wheel-free} if it has no induced $k$-wheel. A graph is
\emph{chordal} if it has no induced cycle of length at least four.
For a positive integer $k$, we simply write $[k]$ to denote the set $\{1,2,\ldots,k\}$, and we say an index $i\in [k]$, if $i\in \{1,2,\ldots,k\}$ and $i$ modulo $k$.

For a set $U\subset V(G)$, let $R_U$ denote a maximum stable set of $U,$ if $U\neq\es$, otherwise let $R_U:=\es$.

\smallskip
Let $G$ be a graph. Suppose $X$ is a subset of $V(G)$ that induces a $P_3$-free graph in $G$.  Then each component of $G[X]$ is a complete subgraph of $G$, and so the set $X$ can be written as a disjoint union  of (nonempty) cliques;
each such clique is a maximal clique of $G[X]$, and  we refer to such cliques as `{\it $X$-cliques}'.
We say that a set $S\subseteq V(G)\sm X$ is \emph{complete to exactly one $X$-clique}, if there is an $X$-clique, $K$, such that $S$ is complete to $K$, and anticomplete to $X\sm K$. Let $v\in V(G)\sm X$ be any vertex.  We say that the vertex $v$ is \emph{good with respect to $X$} if it satisfy the following two conditions:
$(a)$ If $v$ has a neighbor in an $X$-clique, say $K$, then $v$ is complete to $K$, and $(b)$ $v$ is complete to at least one $X$-clique.

\smallskip
We use the following   simple observations often.

\vspace{-0.25cm}
\newcounter{ici}
\begin{enumerate}[label=(O\arabic*),series=edu*]\itemsep=0pt
\item\label{p5} Let $G$ be a $P_5$-free graph. Let $A$, $B_1$ and $B_2$ be three disjoint, nonempty, and mutually anticomplete subsets of $V(G)$.  Let $x$ and $y$ be two nonadjacent vertices in $V(G)\sm (A\cup B_1\cup B_2)$ such that  $x$ and $y$ have a common neighbor in $A$, $x$ has a neighbor in $B_1$, and $y$ has a neighbor in $B_2$. Then $x$ and $y$ must have a common neighbor in either $B_1$ or $B_2$.

\item  \label{P3free}
 Let $G$ be a $4$-wheel-free graph. Let $S$ be a subset of $V(G)$. If there are nonadjacent vertices $u,v\in V(G)\sm S$ such that $\{u,v\}$ is complete to $S$, then $S$ induces a $P_3$-free graph.

\item \label{new-obs}
Let $G$ be any graph. Let $D_1$, $D_2$ and $D_3$ be three disjoint, nonempty subsets of $V(G)$, each induces a $P_3$-free graph.  Suppose that each $D_i$-clique is either complete or anticomplete to a $D_j$-clique, where $i\neq j\in \{1,2,3\}$. If $M$ is a maximal clique in $G$ containing at least one vertex from each $D_1$ and $D_2$, then $R_{D_1}\cup R_{D_2}$ meets $M$ twice.

{\it Proof of \ref{new-obs}}.~If $M\cap D_3=\es$, then clearly the assertion holds. So we may assume that $M\cap D_3\neq \es$.
Then by our assumption, $M$ is of the form $\cup_{i=1}^3D_i^*$, where $D_i^*$ is a $D_i$-clique. Since $R_{D_1}$ contains a vertex from $D_1^*$, and $R_{D_2}$ contains a vertex from $D_2^*$,  we conclude that $R_{D_1}\cup R_{D_2}$ meets $M$ twice. This proves \ref{new-obs}. $\Diamond$

\end{enumerate}

 Next we prove a structure theorem  for a subclass of ($P_5$,\,$4$-wheel)-free graphs, namely   the class of (\,$3K_1$,\,$4$-wheel)-free graphs,   and use it later.

\begin{lemma} \label{3K1W4-free-Struc}
If $G$ is a (\,$3K_1$,\,$4$-wheel)-free graph, then $G$ is either a quasi-line graph or a nice graph.
\end{lemma}
\no{\it Proof of Lemma~\ref{3K1W4-free-Struc}}.~Let $G$ be a ($3K_1$,\,$4$-wheel)-free graph, and let $v\in V(G)$ be arbitrary. First suppose that $G[N(v)]$ is chordal.
 Since the complement graph of a $3K_1$-free chordal graph is a ($K_3, 2K_2, C_5$)-free graph (which is a bipartite graph), we see that $N(v)$ can be expressed as  union of two cliques, and hence $G$ is a quasi-line graph,  since $v$ is arbitrary. So we may assume that  $G[N(v)]$ is not chordal. Then since $G$ does not contain a $4$-wheel, $G[N(v)]$ contains an induced $C_k$ for some $k\geq 5$. Since, for $k\geq 6$,  $C_k$ contains an induced $3K_1$, $G[N(v)]$ contains an induced $C_5$, say  $C$. Hence  $G$ contains a $5$-wheel, induced by the vertices $V(C)\cup \{v\}$. Then it is shown in  Theorem~3 (Case 1.1) of \cite{CKS2008}  that $G$ is a clique-blowup of a $5$-wheel, say  $W$ with the vertex-set $\{v_1, v_2,\ldots, v_6\}$ such that $v_1$-$v_2$-$\cdots$-$v_5$-$v_1$ is a $C_5$, and $v_6$ is complete to $\{v_1, \ldots,v_5\}$. By the definition of clique-blowup, $V(G)$ is partitioned into cliques $Q_{v_i}$, $v_i \in V(W)$.
Then clearly $S_1:=R_{Q_{v_1}} \cup R_{Q_{v_3}}$, $S_2:=R_{Q_{v_2}} \cup R_{Q_{v_4}}$, and $S_3:=R_{Q_{v_5}}$ are three stable sets in $G$ such that $\omega(G-(S_1\cup S_2\cup S_3)) \leq \omega(G)-2$, and so $G$ is nice. This proves Lemma~\ref{3K1W4-free-Struc}. \hfill{$\Box$}

\smallskip
We will  also use the following lemma.
\begin{lemma}[\cite{KM-SIDMA}]\label{lem:squarefree}	
	Let $G$ be a graph. Let $A$ and  $B$ be two disjoint cliques such that $G[A\cup B]$ is $C_4$-free. If every vertex in $A$ has a neighbor in $B$, then some vertex in $B$ is complete to $A$.
\end{lemma}

For the reader's convenience, we give a sketch of the graph (in most cases), and we use the following representations:
The shapes (circles or ovals) represent a collection of sets into which the vertex-set of the graph is partitioned. The sets inside an oval form a partition of that set.
Each shaded shape represents a nonempty clique, and other shapes induce a $P_3$-free subgraph. A solid line between any two shapes represents that the respective sets are complete to each other. A dashed line between any two shapes represents that the adjacency between these sets are arbitrary, but are restricted with some conditions. The absence of a line between any two shapes represents that the respective sets are anticomplete to each other.

\section{Structural properties of ($P_5$,\,$4$-wheel)-free atoms with a $C_5$}\label{sec-C5}
In this section, we present some important and useful structural properties of  ($P_5$, $4$-wheel)-free atoms with an induced $C_5$, and use them in Section~\ref{sec-mainstruc}.

 Let $G$ be a connected ($P_5$, $4$-wheel)-free atom.
Suppose that $G$ contains an induced  $C_5$, say $v_1$-$v_2$-$v_3$-$v_4$-$v_5$-$v_1$.  Then we may assume that there are five nonempty and pairwise disjoint sets $A_1, A_2, \ldots, A_5$ such that for each $i$ modulo~$5$ the
set $A_i$ is complete to $A_{i-1}\cup A_{i+1}$, and anticomplete to
$A_{i-2}\cup A_{i+2}$.  Let $A:=A_1\cup\cdots\cup A_5$.  We choose
these sets such that $A$ is maximal, and let $v_i\in A_i$.  From now on, in this section, every subscript is
understood modulo~$5$.  Let $T:=\{x\in V(G)\setminus A \mid x$ has no
neighbor in $A\}$,  $Z:=\{x\in V(G)\setminus A \mid x$  has a neighbor in each $A_{i}$,  $i\in [5]\}$, and for each $i\in [5]$: let $X_i := \{x\in V(G)\setminus A \mid x \mbox{ has a neighbor in each } A_j, j\in \{i, i+2,i-2\},\mbox{and anticomplete to } A_{i-1}\cup A_{i+1}\}$,   and $Y_i := \{x\in V(G)\setminus A \mid x \mbox{ has a neighbor in each } A_{j}, j\in[5], j\neq i, \mbox{ and anticomplete to } A_{i}\}$.
Let $X:= X_1\cup\cdots \cup X_5$ and $Y:= Y_1\cup\cdots\cup Y_5$. Then:

\begin{prop}\label{part}
  \label{part} $V(G) = A\cup X\cup Y\cup Z\cup T$.
\end{prop}
Suppose there is a vertex $p\in V(G) \setminus (A\cup X\cup Y\cup Z\cup T)$. Since $p\notin T$, $p$ has a neighbor in $A$.  Then since $p\notin X\cup Y\cup Z$, up to symmetry, we have two cases. $(1)$ Suppose $p$ has a neighbor $a_i\in A_i$, and anticomplete to $A_{i-1}\cup A_{i-2}\cup A_{i+2}$. Then $p$-$a_i$-$v_{i-1}$-$v_{i-2}$-$v_{i+2}$ is a $P_5$, a contradiction. $(2)$ Suppose   $p$ has   neighbors   $a_{i-1}\in A_{i-1}$   and    $a_{i+1} \in A_{i+1}$, and anticomplete to $A_{i-2}\cup A_{i+2}$.   Then $p$ is complete to $A_{i-1}$ for otherwise for any nonneighbor of $p$ in $A_{i-1}$, say $b_{i-1}$, $b_{i-1}$-$v_{i-2}$-$v_{i+2}$-$a_{i+1}$-$p$ is a $P_5$, a contradiction. Likewise, $p$ is complete to $A_{i+1}$. But then $p$ can be added to $A_i$ contradicting the maximality of $A$.  These contradictions complete the proof of  \ref{part}. \hfill{$\sq$}

\smallskip
  Moreover, the following statements hold, for each $i\in [5]$:

\begin{prop}

 \label{Ai-P3free}
 $G[A_i]$ is $P_3$-free. So $G[A]$ is a blowup of $C_5$.
\end{prop}
 If $G[A_i]$ contains an induced $P_3$, say $u_1$-$u_2$-$u_3$, then $\{u_1,v_{i+1},u_3,v_{i-1},u_2\}$ induces a $4$-wheel, a contradiction. This proves \ref{Ai-P3free}. \hfill{$\sq$}

\begin{prop}
 \label{XiAi}
 $X_i$ is complete to $A_i$.
\end{prop}
  If there are nonadjacent vertices, say $x\in X_i$ and $p\in A_i$, then for any  neighbor of $x$ in $A_{i+2}$, say $a_{i+2}$, we see that $v_{i-1}$-$p$-$v_{i+1}$-$a_{i+2}$-$x$ is a $P_5$, a contradiction. This proves \ref{XiAi}. \hfill{$\sq$}

\begin{prop}
 \label{Ai+2-Ai-2-K-xi}
 If $K$ is an $A_{i+2}$-clique (or an $A_{i-2}$-clique), then any $x\in X_i$ which has a neighbor in $K$, is complete to $K$.
In particular, if $A_{i+2}$ is a clique, then  $X_i$ is complete to $A_{i+2}$. Likewise,  if $A_{i-2}$ is a clique, then    $X_i$ is complete to $A_{i-2}$.
\end{prop}
By symmetry, it is enough to prove the assertion for $A_{i+2}$.
If $x\in X_i$ is not complete to $K$, then by assumption, there are vertices $a,b$ in $K$ such that $ab,ax\in E$ and $bx\notin E$. But then by \ref{XiAi}, $b$-$a$-$x$-$v_i$-$v_{i-1}$ is a $P_5$, a contradiction. This proves   \ref{Ai+2-Ai-2-K-xi}. \hfill{$\sq$}

\begin{prop}
  \label{x-good}  Each vertex $X_i$ is good with respect to $A_{i+2}$, and $A_{i-2}$.
  \end{prop}
    The proof of \ref{x-good} follows from the definition of $X_i$, \ref{Ai-P3free} and \ref{Ai+2-Ai-2-K-xi}. \hfill{$\sq$}

\begin{prop}
 \label{xAi+2Ai-2}
 Each vertex in $X_i$ is complete to either $A_{i+2}$ or $A_{i-2}$.
\end{prop}
Let $x\in X_i$, and suppose that the assertion is not true.  Then there are vertices $p\in A_{i+2}$ and $q\in A_{i-2}$ such that $x$ is anticomplete to $\{p,q\}$. By the definition of $X_i$, $x$ has a neighbor in $A_{i+2}$, say $r$. Then by \ref{Ai-P3free} and \ref{Ai+2-Ai-2-K-xi}, $pr\notin E$. But then by \ref{XiAi}, $p$-$q$-$r$-$x$-$v_i$ is a $P_5$, a contradiction. This proves \ref{xAi+2Ai-2}. \hfill{$\sq$}

  \begin{prop}
  \label{nadj-x-cn} Any two nonadjacent vertices in $X_i$ have a common neighbor in $A_{i+2}$, and in $A_{i-2}$.
  \end{prop}
  The proof of \ref{nadj-x-cn} follows from the definition of $X_i$, \ref{Ai-P3free} to \ref{Ai+2-Ai-2-K-xi}, and by \ref{p5}. \hfill{$\sq$}

\begin{prop}\label{2nax-Ai+2Ai-2clq}
If $X_i$ has two nonadjacent vertices which are complete to
$A_{i+2}\cup A_{i-2}$, then $A_{i+2}\cup A_{i-2}$ is a clique.
\end{prop}
Suppose there are nonadjacent vertices in $A_{i+2}$, say $a$ and $a'$.  Let $x,x'$ be two nonadjacent vertices in $X_i$ which are complete to $A_{i+2}\cup A_{i-2}$. Then for any $a''\in A_{i-2}$, $\{x,a,x',a',a''\}$ induces a $4$-wheel, a contradiction. So $A_{i+2}$ is a clique. Likewise, $A_{i-2}$ is a clique. This proves  \ref{2nax-Ai+2Ai-2clq}. \hfill{$\sq$}

\begin{prop}
 \label{xi-pure}
 If some $x\in X_i$ has a neighbor in $T$, then $x$ is complete to $A_{i-2}\cup A_{i+2}$.
\end{prop}
Let   $t\in T$ be a neighbor of $x$. By \ref{XiAi} and \ref{xAi+2Ai-2}, we may assume that $x$ is complete to $A_i\cup A_{i-2}$. If $x$ has a nonneighbor in $A_{i+2}$, say $p$, then $p$-$v_{i+1}$-$v_{i}$-$x$-$t$ is a $P_5$, a contradiction. So $x$ is complete to $A_{i+2}$. This proves   \ref{xi-pure}. \hfill{$\sq$}

\begin{prop}
 \label{Xi-P3free}
 $G[X_i]$ is $P_3$-free.
\end{prop}
 Suppose to the contrary that $G[X_i]$ induces a $P_3$ with vertex-set $\{a_1,a_2,a_3\}$. Then by \ref{xAi+2Ai-2} and by the pigeonhole principle, we may assume that $\{a_1,a_2\}$ is complete to $A_{i-2}$. Also by the definition of $X_i$, $a_3$ has a neighbor in $A_{i-2}$, say $p$. Then by \ref{XiAi}, $\{v_i,a_1,a_2,a_3,p\}$ induces a $4$-wheel, a contradiction. This proves \ref{Xi-P3free}. \hfill{$\sq$}

\begin{prop}
 \label{XiXi+1}
 $X_i$ is complete to $X_{i+1} \cup X_{i-1}$.
\end{prop}
Let $x\in X_{i}$ and $x'\in X_{i+1}$, and suppose that $x,x'$ are nonadjacent. By definition,   pick a neighbor of $x'$ in $A_{i-1}$, say $p$, and a neighbor of $x$ in $A_{i+2}$,  say $q$. Then by \ref{XiAi}, $p$-$x'$-$v_{i+1}$-$q$-$x$ is a $P_5$,  a contradiction. So $X_i$ is complete to $X_{i+1}$. Likewise, $X_i$ is complete to $X_{i-1}$. This proves \ref{XiXi+1}. \hfill{$\sq$}

\begin{prop}\label{XiAi-2clqcor} {\it Let $K$ be an $X_i$-clique. Then the following hold: \begin{enumerate}[label=(\alph*)]\itemsep=0pt
 \item  Suppose that there is a vertex $x\in X_{i+2}$ which is anticomplete to $K$, and $Q$ is an $A_{i-2}$-clique such that $N(K)\cap Q\neq \es$. Then $K$ is complete to $Q$.
     \item Suppose that there is a vertex $x\in X_{i-2}$ which is anticomplete to $K$, and  $Q'$ is an $A_{i+2}$-clique such that $N(K)\cap Q'\neq \es$. Then $K$ is complete to $Q'$.\end{enumerate}}
\end{prop}
   We prove $(a)$, and the proof of $(b)$ is similar. Suppose that the assertion is not true. Then there are vertices $p\in K$ and $r\in Q$ such that $pr\notin E$. By assumption, there is a vertex $q\in K$ such that $q$ has a neighbor in $Q$, and hence by \ref{Ai+2-Ai-2-K-xi}, $qr\in E$. Then for any neighbor of $x$ in  $A_{i-1}$, say $a$, we see that $p$-$q$-$r$-$a$-$x$ is a $P_5$, a contradiction. This proves \ref{XiAi-2clqcor}. \hfill{$\sq$}

\begin{prop}\label{KcompK'}{\it   Suppose $K$ is an $X_i$-clique and $K'$ is an $X_{i+2}$-clique such that $K$ is  complete to $K'$. Then the following hold:
 \begin{enumerate}[label=(\alph*)]
 \item $K$ is anticomplete to $X_{i+2}\sm K'$ (likewise, $K'$ is  anticomplete to $X_{i}\sm K$), and  $X_i\sm K$ is anticomplete to $X_{i+2}\sm K'$.

 \item  $K$ is complete to exactly one $A_{i+2}$-clique. Likewise, $K'$ is complete to exactly one $A_i$-clique.

\item $K$ is anticomplete to $X_{i-2}$. Likewise, $K'$ is anticomplete to $X_{i-1}$.
\end{enumerate}}
\end{prop}
$(a)$: Suppose to the contrary that $K$ is not anticomplete to $X_{i+2}\sm K'$. Then there are vertices $u\in K$, $v\in K'$ and $w\in X_{i+2}\sm K'$ such that $uv,uw\in E$ and $vw\notin E$.  Then by \ref{nadj-x-cn}, $v$ and $w$ have a common neighbor in $A_i$, say $p$. But then for any neighbor of $u$  in $A_{i+2}$, say $q$, by \ref{XiAi}, $\{p,v,q,w,u\}$ induces a $4$-wheel, a contradiction. So $K$ is anticomplete to $X_{i+2}\sm K'$. Likewise, $K'$ is  anticomplete to $X_{i}\sm K$. This proves the first assertion of $(a)$.

To prove the second assertion of $(a)$, suppose there are adjacent vertices, say $u'\in X_i\sm K$ and $v'\in X_{i+2}\sm K'$. Then for any $v\in K'$, since $vv'\notin E$, by \ref{nadj-x-cn}, $v$ and $v'$ have a common neighbor in $A_{i-1}$, say $p$. But then for any $u\in K$, by using the first assertion of $(a)$, we see that $u$-$v$-$p$-$v'$-$u'$ is a $P_5$, a contradiction. This proves the second assertion of $(a)$.

\no $(b)$: First we show that each vertex in $K$ is complete to exactly one $A_{i+2}$-clique. Suppose not. Then by \ref{x-good}, there are  vertices $p\in K$ and $a, a'\in A_{i+2}$ such that $pa,pa'\in E$ and $aa'\notin E$. But then for any $q\in K'$, and for any neighbor of $p$ in $A_{i-2}$, say $r$, by \ref{XiAi}, $\{r,a,q,a',p\}$ induces a $4$-wheel, a contradiction. So each vertex in $K$ is complete to exactly one $A_{i+2}$-clique.  Now we show that $K$ is complete to exactly one $A_{i+2}$-clique. Suppose not. Then by \ref{Ai+2-Ai-2-K-xi} and by the earlier argument, there are vertices $u,v\in K$ and $p \in A_{i+2}$ such that $up\in E$ and $vp\notin E$. Then by \ref{xAi+2Ai-2}, $v$ is complete to $A_{i-2}$. But then for any neighbor of $u$ in $A_{i-2}$, say $a$, and for any $q\in K'$, by \ref{XiAi}, $\{a,v,q,p,u\}$ induces a $4$-wheel, contradiction. This proves $(b)$.

\no $(c)$: Let $u\in K$ and $v\in X_{i-2}$, and suppose $u,v$ are adjacent. Let $r\in K'$. By \ref{XiXi+1}, $v$ and $r$ are adjacent. Now pick any neighbor of $u$ in $A_{i+2}$, say $p$, and in  $A_{i-2}$, say $q$. Then by \ref{XiAi}, $\{p,q,v,r,u\}$ induces a $4$-wheel, a contradiction. This proves $(c)$.

This completes the proof of \ref{KcompK'}. \hfill{$\sq$}

\begin{prop}\label{G7-clm1}
Let $K$ be an $X_i$-clique and  $K'$ be an $X_{{i-1}}$-clique. If $Q$ is an $A_{i+2}$-clique such that $N(K)\cap Q\neq \es$ and $N(K')\cap Q\neq \es$, then $K\cup K'$ is complete to $Q$.
\end{prop}
We prove the assertion for $i=1$. Suppose that $K$ is not complete to $Q$. Then there are  vertices $p\in K$ and $r\in Q$ such that $pr\notin E$.  By assumption, there is a vertex $q\in K$ such that $q$ has a neighbor in $Q$, and so by \ref{Ai+2-Ai-2-K-xi},  $qr\in E$.  Also by our assumption, there is a vertex $w\in K'$ such that $w$ has a neighbor in $Q$, and again by \ref{Ai+2-Ai-2-K-xi}, $wr\in E$. Since $p$ is not complete to $A_{3}$,  $p$ is complete to $A_{4}$, and so  $p$ and $q$ share a common neighbor in $A_{4}$, say $x$. Then since $X_1$ is complete to $X_5$ (by \ref{XiXi+1}), we see that $\{w,r,x,p,q\}$ induces a $4$-wheel,  a contradiction. So $K$ is complete to $Q$. Likewise, $K'$ is complete to $Q$. This proves \ref{G7-clm1}.  \hfill{$\sq$}

\medskip

For each $i\in [5]$, if $X_i\neq \es$, let  ${\mathbb{W}}_i$ denote the set $\{X^*\cup A^* \mid X^*$  is an $X_i$-clique and $A^*$  is an $A_i$-clique  such that $|X^*\cup A^* |= \omega(G)\}$, otherwise let ${\mathbb{W}}_i:=\es$. Next we have the following:

\begin{prop}\label{3-XiXi+1-Ai-2}
Let $K$ be an $X_i$-clique and $K'$ be an $X_{i+1}$-clique, and let $A_i^*$ be an $A_i$-clique, and  $A_{i+1}^*$ be an $A_{i+1}$-clique. Suppose $ K\cup A_i^*\in {\mathbb{W}}_i$ and $K'\cup A_{i+1}^*\in {\mathbb{W}}_{i+1}$. Then for any $A_{i+2}$-clique $D_{i-2}$,  $K\cup K'\cup D_{i-2}$ is not a clique.
\end{prop}
By \ref{XiXi+1}, $K\cup K'$ is a clique. Suppose  there is an $A_{i-2}$-clique, say $D$, such that $K\cup K'\cup D$ is a clique.  Let $q:=\omega(G)$. Then $|K\cup K'|< q$ (since $D\neq \es$). Then since  $|K\cup A_i^*|+|K'\cup A_{i+1}^*|=2q$, we have $2q=|A_i^*\cup A_{i+1}^*|+|K\cup K'|<|A_i^*\cup A_{i+1}^*|+q$, and hence $|A_i^*\cup A_{i+1}^*|> q$,  a contradiction since $A_i^*\cup A_{i+1}^*$ is a clique. This proves \ref{3-XiXi+1-Ai-2}. \hfill{$\sq$}

\begin{prop}
\label{xi-compT}
 Let $Q$ be the vertex-set of a component of $G[T]$. Then each vertex in $X_i$ is either complete  or anticomplete to $Q$.
\end{prop}
Otherwise,  there are adjacent vertices $q, q'$ in $Q$, and a vertex $x\in X_i$ such that $xq\in E$ and $xq'\notin E$; but then by \ref{XiAi}, $q'$-$q$-$x$-$v_i$-$v_{i-1}$ is a $P_5$, a contradiction. This proves \ref{xi-compT}. \hfill{$\sq$}

\begin{prop}
 \label{SA-Yi-Ai-1-Ai+1}
 For $j\in \{i-1,i+1\}$, if $A_j$ is not a clique, then $Y_i$ is complete to $A_j$.
\end{prop}
 We may assume, up to symmetry,  that $j=i+1$. Let $y\in Y_i$. Then by the definition of $Y_i$,
 $y$ has a neighbor in $A_{i+1}$, say $p$. Let $K$ be the $A_{i+1}$-clique containing $p$. Since $A_{i+1}$ is not a clique,
 $A_{i+1}\sm K \neq \es$. Now if $y$ is nonadjacent to some $q\in A_{i+1}\sm K$ (say),  then for any neighbor of $y$  in $A_{i-2}$, say $r$,  we see that $q$-$v_i$-$p$-$y$-$r$ is   a $P_5$, a contradiction; so $y$ is complete to $A_{i+1}\sm K$. By the same argument, since $A_{i+1}\sm K$ is nonempty, $y$ is complete to $K$. This proves \ref{SA-Yi-Ai-1-Ai+1}, since $y$ is arbitrary. \hfill{$\sq$}

\begin{prop}
 \label{yi-Ai-1-Ai+1}
 Each vertex in $Y_i$ is  complete to either $A_{i-1}$ or $A_{i+1}$.
\end{prop}
Let $y\in Y_i$ . Suppose  $y$ has a nonneighbor in each $A_{i-1}$ and $A_{i+1}$, say $a$ and $a'$ respectively.  So by \ref{SA-Yi-Ai-1-Ai+1}, $A_{i-1}$ and $A_{i+1}$ are cliques. Now by the definition of $Y_i$, pick any neighbor of $y$ in each $A_{i-1}$ and $A_{i+1}$, say $b$ and $b'$ respectively. Then  $a$-$b$-$y$-$b'$-$a'$ is a $P_5$, a contradiction. This proves \ref{yi-Ai-1-Ai+1}. \hfill{$\sq$}

\begin{prop}\label{yi-compT}
Let $Q$ be the vertex-set of a component of $G[T]$. Then each vertex in $Y_i$ is either complete  or anticomplete to $Q$.
\end{prop}
Otherwise,  there are adjacent vertices $q, q'$ in $Q$, and a vertex $y\in Y_i$ such that $yq\in E$ and $yq'\notin E$; but then for any neighbor of $y$ in $A_{i+1}$, say $a$, we see that $q'$-$q$-$y$-$a$-$v_i$ is a $P_5$, a contradiction. This proves \ref{yi-compT}. \hfill{$\sq$}

\begin{prop} \label{5WF-T-P3free}If $Z=\es$, then $G[T]$ is $P_3$-free. \end{prop}
 Suppose  that there is a component of $G[T]$ which has an induced $P_3$, say $t_1$-$t_2$-$t_3$, and let $Q$ be the vertex-set of that component.  Since $G$ is connected, and since $N(Q)\cap (X\cup Y)$ is not a clique cutset,
there are nonadjacent vertices in $N(Q)\cap (X\cup Y)$, say $u$ and $v$. Then by \ref{xi-compT} and \ref{yi-compT}, $\{u,v\}$ is complete to $Q$; but then $\{u,t_1,v,t_3,t_2\}$ induces a $4$-wheel, a contradiction. This proves \ref{5WF-T-P3free}. \hfill{$\Box$}

\smallskip
So if $Z=\es$, then, by \ref{5WF-T-P3free}, \ref{xi-compT} and \ref{yi-compT}, each vertex in $X\cup Y$ is either anticomplete or good with respect  to $T$.

\begin{prop} \label{TXi-nonadj-last}   Suppose there are vertices  $t\in T$, $u\in X_i$ and $v\in X_{i-2}\cup X_{i+2}\cup Y_i\cup Y_{i+1}\cup Y_{i-1}\cup Z$ such that $ut\in E$  and $uv\notin E$. Let $K$ be the $X_i$-clique containing $u$. Then  the following hold:
\begin{enumerate}\itemsep=0pt
    \item[(a)]  $t$ is adjacent to $v$.
    \item[(b)] If $v$ is  anticomplete to $K$,  then $t$ is complete to $K$. Moreover, if $T^*$ is the component of $T$ containing $t$, then $T^*$ is complete to $K$.
\end{enumerate}
\end{prop}
First note that  $v$ has a neighbor  in one of $A_{i-1}$, $A_{i+1}$. We may assume, up to symmetry, that $v$ has a neighbor in $A_{i-1}$, say $p$. So $v\notin X_{i-2}\cup Y_{i-1}$.

 \no $(a)$: Suppose $t$ is nonadjacent to $v$. If $v$ is nonadjacent to some vertex in $A_i$, say $q$, then, by \ref{XiAi}, $v$-$p$-$q$-$u$-$t$ is a $P_5$, a contradiction; so $v$ is complete to $A_i$. Thus  $v\notin Y_i$, and so  $v\in X_{i+2} \cup Y_{i+1}\cup Z$. Then since $ut\in E$, by \ref{xi-pure}, $u$ is complete to $A_{i+2}$, and so $u$ and $v$ have a common neighbor in $A_{i+2}$, say $r$. But then $t$-$u$-$r$-$v$-$p$ is a $P_5$, a contradiction. This proves $(a)$.

\no $(b)$: If there is a vertex $u'\in K$ such that $u't\notin E$, then by $(a)$, $u'$-$u$-$t$-$v$-$p$ is a $P_5$, a contradiction. This proves the first assertion of $(b)$. The second assertion of $(b)$ follows from \ref{xi-compT}. \hfill{$\sq$}

\section{Structure of ($P_5$,\,$4$-wheel)-free atoms}\label{sec-mainstruc}
In this section, we give a proof of Theorem~\ref{thm:struc-new}. Indeed, we prove each of the item in Theorem~\ref{thm:struc-new} separately, and are given below in Theorems~\ref{thm-5wheel}, \ref{thm:P5wheel-free-C5} and \ref{thm-C7c}  respectively. In most cases we show that $G$ is nice, and to do the same it is enough to find three stable sets  $S_1, S_2$, and $S_3$ such that $S_1\cup S_2\cup S_3$ meets each maximum clique of $G$ at least twice, and other maximal cliques at least once.

\subsection{Structure of ($P_5$,\,$4$-wheel)-free atoms with an induced $5$-wheel}

Let $G$ be a connected ($P_5$,\,$4$-wheel)-free atom that contains an induced $5$-wheel, say with the $5$-cycle $v_1$-$v_2$-$v_3$-$v_4$-$v_5$-$v_1$ plus a vertex $z^*$ that is adjacent to $v_i$, for all $i\in [5]$.     Then we define the sets $A$, $X$, $Y$, $Z$ and $T$ as in Section~\ref{sec-C5} with $v_i\in A_i$ for each $i\in [5]$, and  we use the facts \ref{part}--\ref{TXi-nonadj-last} shown in Section~\ref{sec-C5}.  Note that $z^*\in Z$.  Moreover, the graph $G$ has some more structural properties, and are given in \cref{prop-Z,prop-ZXY,prop-ZT} below.

\begin{lemma}\label{prop-Z}
For $i\in [5]$, the following properties hold:

\begin{lemmalist}
\item \label{z-Ai-K} Let $K$ be an $A_i$-clique. If a vertex in $Z$ has a neighbor in $K,$ then  it is  complete to $K,$ and anticomplete to $A_i\sm K$.
More precisely, each vertex in $Z$ is complete to exactly one $A_i$-clique.
\item  \label{Z-AjAj+2Aj-2}
There is an index $j\in [5]$ such that $A_j$, $A_{j-2}$ and $A_{j+2}$ are cliques.
\item   \label{Zclq}
 $Z$ is a clique.
\item\label{Ai*-Z}
 There is an $A_i$-clique, say  $A_i^*$, such that  $Z$ is complete to $A_i^*$, and anticomplete to $A_i\sm A_i^*$.
\item\label{Z-O3}
 $R_{A_{i+2}}\cup R_{A_{i-2}}$ meets each maximal clique of $G$ in   $G[Z\cup A_{i+2}\cup A_{i-2}]$ twice.

\end{lemmalist}
\end{lemma}
\no{\it Proof of \cref{prop-Z}}.~$(i)$:~Let $z\in Z$, and suppose $z$ has a neighbor in $K$, say $p$. Pick a neighbor of $z$ in each $A_{i+1}$ and $A_{i-1}$, say $a$ and $a'$ respectively.
If there is a vertex $q\in K$ which is nonadjacent to $z$, then since $K$ is a clique, $pq\in E$, and then $\{a,z,a',q,p\}$ induces a $4$-wheel, a contradiction; so $z$ is complete to $K$. Next, if there is a vertex $r\in A_i\sm K$ which is adjacent to $z$, then clearly $pq\notin E$, and then $\{p,a,q,a',z\}$ induces a $4$-wheel, a contradiction; So $z$ is anticomplete to $A_i\sm K$. This proves \cref{z-Ai-K}.  $\Diamond$

\smallskip
\no{$(ii)$}:~We first show that, for $i\in [5]$, each vertex in $Z$  is complete to either  $A_i$ or $A_{i+1}$. Suppose   not. Then there are vertices $b\in A_i$ and $b'\in A_{i+1}$ such that $zb,zb'\notin E$. Now pick a neighbor of    $z$ in each $A_i$ and $A_{i-2}$, say $a$ and $a'$, respectively. Then by \cref{z-Ai-K},  $ab\notin E$; but then $b$-$b'$-$a$-$z$-$a'$ is a $P_5$, a contradiction. So  each vertex in $Z$  is complete to either  $A_i$ or $A_{i+1}$. Then for $i\in [5]$,  since $z^*\in Z$  is complete to exactly one $A_i$-clique (by \cref{z-Ai-K}), we see that  either $A_i$ or $A_{i+1}$ is a clique, and so \cref{Z-AjAj+2Aj-2} holds.  $\Diamond$

\smallskip
\no{$(iii)$}:~Suppose there are nonadjacent vertices, say $z_1,z_2$ in $Z$. Then by \cref{Z-AjAj+2Aj-2}, we may assume that $A_1$ and $A_{3}$ are cliques. So by \cref{z-Ai-K}, $\{z_1,z_2\}$ is complete to $A_1\cup A_3$. Also, by the definition of $Z$,  \ref{p5} and  by \cref{z-Ai-K},  it follows that $z_1$ and $z_2$ have a common neighbor in $A_2$, say $p$. Then $\{v_1,z_1,v_3,z_2,p\}$ induces a $4$-wheel, a contradiction. This proves \cref{Zclq}.  $\Diamond$

\smallskip
\no{$(iv)$}:~By \cref{Z-AjAj+2Aj-2}, we may assume that $A_1$, $A_3$ and $A_4$ are cliques. So by \cref{z-Ai-K}, for $j\in \{1,3,4\}$, $A_j$ is our required $A_j^*$.
  This implies that $Z$ is complete to $A_j$, for $j\in \{1,3,4\}$.
  Next we prove that $A_2^*$ and $A_5^*$ exist. Suppose, up to symmetry,   $A_2^*$ does not exist. Then by \cref{z-Ai-K},    there are vertices $z_1,z_2\in Z$, and a vertex  $p\in A_2$ such that $z_1p\in E$ and $z_2p\notin E$. By \cref{Zclq}, $z_1z_2\in E$. Then $\{v_1,p,v_3,z_2,z_1\}$ induces a $4$-wheel, a contradiction. So $A_2^*$  exists.  This proves \cref{Ai*-Z}. $\Diamond$

\smallskip
\no{$(v)$}:~The proof  follows from \ref{Ai-P3free}, \cref{Zclq} and \cref{Ai*-Z}. \hfill{$\Box$}

\medskip
Throughout this subsection, for $i\in [5]$,   $A_i^*$ is an $A_i$-clique as in \cref{Ai*-Z}. Note that by \cref{Ai*-Z}, since $z^*\in Z$,   for  $i\in [5]$, we have $v_i\in A_i^*$,    $Z$ is complete to $\{v_1,v_2,\ldots, v_5\}$, and if $A_i$ is a clique, then $A_i=A_i^*$ and $Z$ is complete to $A_i$.

\begin{lemma}\label{prop-ZXY}
For $i\in [5]$, the following properties hold:
\begin{lemmalist}
\item\label{Xi-Z}
$X_i$ is anticomplete to $Z$.
\item \label{Z-Xi-Ai-2Ai+2}
For  $j\in \{i-2,i+2\}$, $X_i$ is   complete to $A_j^*$, and   anticomplete to $A_j\sm A_j^*$.
\item \label{XiZ-O3}
$R_{A_{i+2}}\cup R_{A_{i-2}}$ meets each maximum clique of $G$ in   $G[X_i\cup A_{i+2}\cup A_{i-2}]$ twice.
  \item\label{Z-XiXi+2}
$X_i$ is anticomplete to $X_{i+2}\cup X_{i-2}$.
\item\label{Z-Yemp}
$Y$ is  empty.
\item\label{xi-T-Ai+2}
  If a vertex in $X_i$ has a neighbor in $T$, then $A_{i-2}$ and $A_{i+2}$ are cliques.
\end{lemmalist}
\end{lemma}
\no{\it Proof of \cref{prop-ZXY}}.~$(i)$:~Let $x\in X_i$ and $z\in Z$, and suppose $x,z$ are adjacent. By \ref{XiAi} and \ref{xAi+2Ai-2}, we may assume that $x$ is complete to $A_i\cup A_{i+2}$.  Then $\{v_i,v_{i+1},v_{i+2},x,z\}$ induces a $4$-wheel, a contradiction. This proves \cref{Xi-Z}. $\Diamond$

\smallskip
\no{$(ii)$}:~By \cref{Z-AjAj+2Aj-2}, we may assume that $A_{i-2}$ is a clique; so $A_{i-2}=A_{i-2}^*$. Then by \ref{Ai+2-Ai-2-K-xi}, $X_i$ is complete to $A_{i-2}$. Next we prove for $j=i+2$. Pick any $x\in X_i$. Then by \cref{Xi-Z}, $z^*x\notin E$. Also by \ref{XiAi}, $x$ and $z^*$ have a common neighbor in $A_{i}$. Then by definitions of $X_i$ and $Z$, \cref{Ai*-Z}, and by \ref{p5}, $x$ and $z^*$ must have a common neighbor in $A_{i+2}^*$, say  $p$. So  by \ref{Ai+2-Ai-2-K-xi}, $x$ is complete to $A_{i+2}^*$. Next, if $x$ is adjacent to some vertex in $A_{i+2}\sm A_{i+2}^*$, say $q$, then   $q$-$x$-$p$-$z^*$-$v_{i-1}$ is a $P_5$, a contradiction. So $x$ is anticomplete to $A_{i+2}\sm A_{i+2}^*$.  This proves \cref{Z-Xi-Ai-2Ai+2}, since $x\in X_i$ is arbitrary.  $\Diamond$

\smallskip
\no{$(iii)$}:~Since  $R_{A_{i+2}}$ contains a vertex of $A_{i+2}^*$, and  $R_{A_{i-2}}$ contains a vertex of $A_{i-2}^*$, the proof follows from \ref{Ai-P3free}, \ref{Xi-P3free}, and from \cref{Z-Xi-Ai-2Ai+2}. $\Diamond$

\smallskip
\no{$(iv)$}:~Let $x\in X_i$ and $x'\in X_{i+2}$, and suppose $x,x'$ are adjacent. By \cref{Xi-Z}, $z^*$ is nonadjacent to both $x$ and $x'$, and by \cref{Ai*-Z} and \cref{Z-Xi-Ai-2Ai+2},  $v_{i-2}$ is complete to $\{x,z^*\}$. But now  $v_{i+1}$-$z^*$-$v_{i-2}$-$x$-$x'$ is a $P_5$, a contradiction. This proves \cref{Z-XiXi+2}. $\Diamond$

\smallskip
\no{$(v)$}:~Suppose not, and let $y\in Y_i$.  Then by \ref{SA-Yi-Ai-1-Ai+1} and \cref{Ai*-Z}, $y$ and $z^*$ have a common neighbor in both $A_{i+1}$ and  $A_{i-1}$, say $p$ and $q$, respectively. If $z^*y \in E$, then  $\{y,q,v_i,p,z^*\}$ induces a $4$-wheel, a contradiction; so we may assume that $z^*y\notin E$. Then by definitions of $Y_i$ and $Z$, \ref{SA-Yi-Ai-1-Ai+1}, \cref{Ai*-Z} and by \ref{p5}, $z^*$ and $y$ must have a common neighbor in $A_{i+2}^*$, say $b$, and in $A_{i-2}$, say $a$. Then $\{p,y,a,z^*,b\}$ induces a $4$-wheel, a contradiction. This proves \cref{Z-Yemp}. $\Diamond$

\smallskip
\no{$(vi)$}:~Let $x\in X_i$ be a vertex such that $x$ has a neighbor in $T$. Then by \ref{xi-pure}, $x$ is complete to $A_{i-2}\cup A_{i+2}$. Now the conclusion follows from \cref{Z-Xi-Ai-2Ai+2}. \hfill{$\Box$}

  \begin{lemma}\label{prop-ZT}
 The following properties hold:
\begin{lemmalist}
\item\label{ZTX}
Let  $Q$ be the vertex-set of a  component of $G[T]$. Then there is an index $j\in [5]$ such that $N(Q)\cap X_j$ is nonempty, and  is complete to $Q$. In particular, every vertex in $T$ has a neighbor in $X.$
\item\label{ZT}
$Z$ is complete to $T$.
\item\label{ZT-P3free}
$G[T]$ is $P_3$-free.
\end{lemmalist}
\end{lemma}
\no{\it Proof of \cref{prop-ZT}}.~$(i)$:~We know, by \cref{Z-Yemp}, that $Y=\es$. Since $Z$ is a clique (by \cref{Zclq}), and since $N(Q)\cap Z$ is not a clique cutset, we see that $N(Q)\cap X\neq \es$. So there is an index $j\in [5]$ such that $N(Q)\cap X_j\neq \es$.
Pick any $x\in N(Q)\cap X_j$. Then, by \ref{xi-compT},  $x$ is complete to $Q$.  This proves \cref{ZTX}. $\Diamond$

 \smallskip
\no{$(ii)$}:~Since $X$ is anticomplete to $Z$ (by \cref{Xi-Z}), the proof follows from  \cref{ZTX} and \ref{TXi-nonadj-last}:$(a)$.  $\Diamond$

\smallskip
\no{$(iii)$}:~Let $Q$ be the vertex-set of a component of $G[T]$. Suppose to the contrary  that $G[Q]$ contains a $P_3$, say $p$-$q$-$r$.  Then by \cref{ZTX}, for some $j\in [5]$, there is a vertex $x\in X_j$ which is complete to $Q$. But then by \cref{Xi-Z} and \cref{ZT}, $\{p,z^*,r,x,q\}$ induces a $4$-wheel, a contradiction. This proves \cref{ZT-P3free}, since $Q$ is arbitrary. \hfill{$\Box$}

\begin{figure}[h]
\centering
 \includegraphics[height=3.5cm, width=13cm]{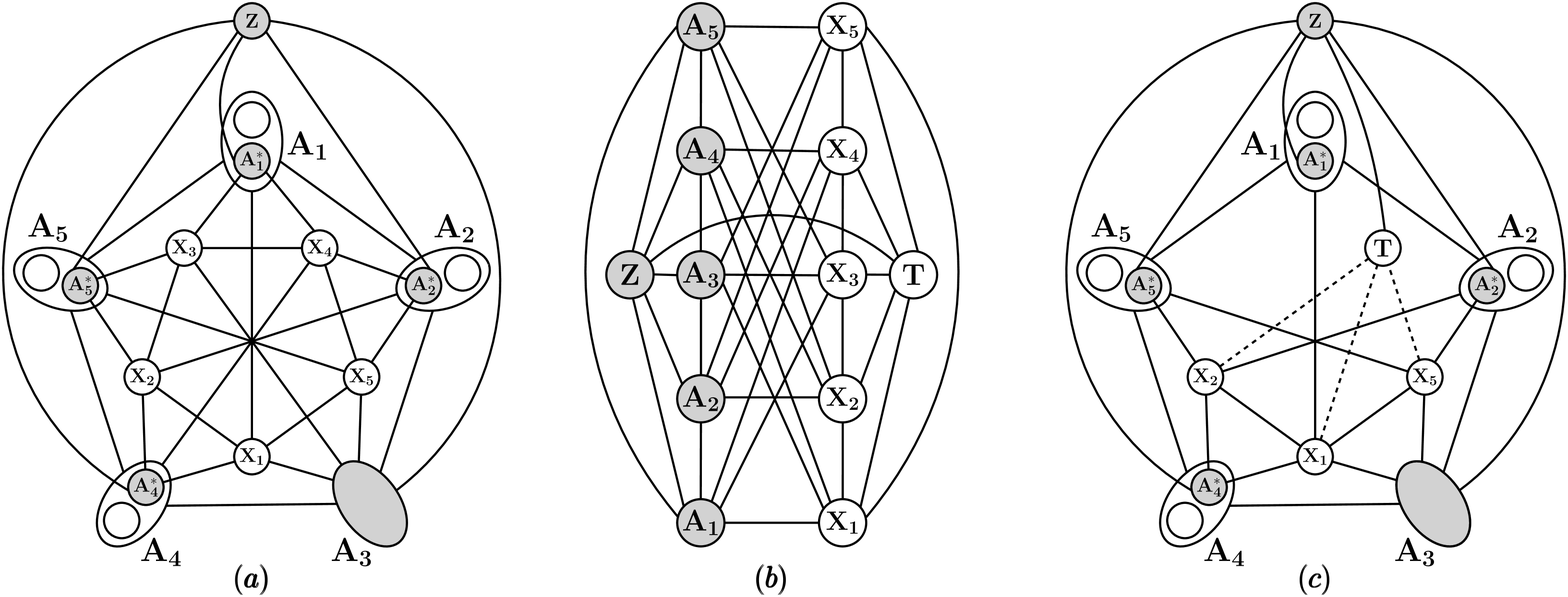}
\caption{Sketch of the graph $G$ in Theorem~\ref{thm-5wheel}: $(a)$ When $T=\es$. $(b)$  Case~1.  (c) Case~2. }\label{fig:Z}
\end{figure}

\smallskip
Now we give our main result of this subsection, and is given below.

\begin{thm}\label{thm-5wheel}
If a connected ($P_5$,\,$4$-wheel)-free  atom $G$ contains an induced $5$-wheel, then $G$ is nice.
\end{thm}
\no{\it Proof of Theorem~\ref{thm-5wheel}}.
Let $G$ be a connected ($P_5$,\,$4$-wheel)-free atom that contains an induced $5$-wheel, say with the $5$-cycle $v_1$-$v_2$-$v_3$-$v_4$-$v_5$-$v_1$ plus a vertex $z^*$ that is adjacent to $v_i$, for all $i\in [5]$.     Then we define the sets $A$, $X$, $Y$, $Z$ and $T$ as in Section~\ref{sec-C5} with $v_i\in A_i$ for each $i\in [5]$. We use the facts \ref{part}--\ref{TXi-nonadj-last} shown in Section~\ref{sec-C5}, and properties in Lemmas~\ref{prop-Z} to \ref{prop-ZT}. Let $\cal M$ denote the set of   maximum cliques in  $G$. To prove the theorem, it is enough to find three stable sets  $S_1, S_2$, and $S_3$ such that $S_1\cup S_2\cup S_3$ meets each maximum clique of $G$ at least twice, and other maximal cliques at least once.
  First suppose that $T=\es$.  (See Figure~\ref{fig:Z}:$(a)$ for a sketch of $G$.) By \cref{Z-AjAj+2Aj-2} and up to relabeling, we may assume that  $A_3$ is a clique.
Then, by \ref{3-XiXi+1-Ai-2}, one of ${\mathbb{W}}_1, {\mathbb{W}}_5$ is empty. We may assume that ${\mathbb{W}}_5=\es$, and so either $X_5=\es$ or  no maximum clique of $G$ is in $G[X_5\cup A_5]$.
  Now we let $S_1:=R_{A_1}\cup R_{A_3}\cup R_{X_2}$, $S_2:=R_{A_2}\cup R_{A_4}\cup R_{X_3}$,   $S_3:=R_{A_5}\cup R_{X_1}\cup R_{X_4}$, and let $S:=S_1\cup S_2\cup S_3$. Clearly $S_1, S_2$ and $S_3$ are stable sets. By \cref{Z-Xi-Ai-2Ai+2} and \cref{XiZ-O3},  $S$ meets each maximal clique of $G$ in $G[A\cup X]$  twice.   Also, by \cref{Z-O3}, $S$ meets each maximal clique of $G$ in $G[A\cup Z]$  twice. So, by \cref{Xi-Z}, we see that $S_1,S_2$ and $S_3$ are the required stable sets. Hence we may assume that $T\neq \es$.  By \cref{ZT-P3free}, we know that $G[T]$ is $P_3$-free.
Let $L$ consist of one vertex from each $T$-clique. Let $L'$ consist of one vertex (which is not in $L$)  from each nontrivial $T$-clique; otherwise we let $L':=\es$. Moreover:

\begin{claim}\label{ZTX'}
 Let  $Q$ be a $T$-clique and let $K$ be an $X_i$-clique. Then $Q$ is either complete or anticomplete to $K$.
 \end{claim}
\no{\it Proof of \cref{ZTX'}}.~Since $Z\neq \es$, the proof   follows from \cref{Xi-Z} and from \ref{TXi-nonadj-last}:$(b)$. $\Diamond$

\smallskip

So any maximal clique containing vertices from both an $X_i$-clique $X_i^*$ and a $T$-clique $T^*$ is $X_i^*\cup T^*$.

\begin{claim}\label{ZXTmax}
If $T^*$ is a $T$-clique such that $Z\cup T^*\in {\cal M}$ or $X_i^*\cup T^*\in{\cal M}$, where $X_i^*$ is an $X_i$-clique and $i\in [5]$,  then $|T^*|\geq 2$.
\end{claim}
\no{\it Proof of \cref{ZXTmax}}.~If $Z\cup T^*\in {\cal M}$, then  since $Z\cup A_1^*\cup A_2^*$ is a clique (by \cref{Zclq} and \cref{Ai*-Z}), we have $|Z\cup T^*|\geq |Z\cup A_1^*\cup A_2^*|$, and thus $|T^*|\geq 2$.
Now if  $X_i^*\cup T^*\in {\cal M}$, then since $X_i^* \cup A_{i+2}\cup A_{i-2}$ is a clique (by \cref{Z-Xi-Ai-2Ai+2} and \cref{xi-T-Ai+2}), we have $|X_i^*\cup T^*|\geq | X_i^* \cup A_{i+2}\cup A_{i-2}|$, and so $|T^*|\geq 2$. This proves  \cref{ZXTmax}.  $\Diamond$

\smallskip
 Now we prove the theorem in two cases:

\medskip
\noindent{\bf Case~1} {\it  Suppose there is an index $j\in [5]$, $X_j,X_{j+2},X_{j-2}\neq \es$. }

\medskip
\no Then by \cref{Z-XiXi+2}, \cref{ZTX} and \ref{TXi-nonadj-last},  we  see that $T$ is complete to $X$. So by \cref{xi-T-Ai+2}, for each $i\in [5]$, $A_i$ is a clique (so $A_i=A_i^*$), and so by \ref{Ai+2-Ai-2-K-xi}, $X_i$ is complete to $A_{i-2}\cup A_{i+2}$. See Figure~\ref{fig:Z}:$(b)$ for a sketch of $G$. First suppose that  for any $T$-clique $T^*$, $Z\cup T^*\notin {\cal M}$. Now if there is an index $i\in [5]$ such that  ${\mathbb{W}}_i\neq \es$, then we let  $k=i$, otherwise we let $k=j$. Then since for each $i\in [5]$, $A_i$ is a clique, by \ref{3-XiXi+1-Ai-2}, ${\mathbb{W}}_{k-1}, {\mathbb{W}}_{k+1}=\es$.
So we let $S_1:=R_{A_k}\cup R_{A_{k+2}}\cup R_T$, $S_2:= R_{A_{k+1}}\cup R_{A_{k-2}}\cup R_{X_{k+2}}$, and $S_3:= R_{A_{k-1}}\cup R_{X_{k}}\cup R_{X_{k-2}}$. Then, since $R_T\cup R_{X_k} \cup R_{X_{k+2}}\cup R_{X_{k-2}}$ meets each maximum clique of $G$ in   $G[X\cup Z\cup T]$ twice, and the other maximal cliques in $G[X\cup Z\cup T]$ once,  by \cref{XiZ-O3} and \cref{Z-O3}, we see that $S_1,S_2$ and $S_3$ are  the desired stable sets.
So suppose that there is a $T$-clique  $T^*$ such that $Z\cup T^*\in {\cal M}$. Then, by \cref{ZXTmax}, $|T^*|\geq 2$. Now for any $X_i$-clique $X_i^*$, and for any $X_{i+1}$-clique $X_{i+1}^*$, by \ref{XiXi+1},  $|Z\cup T^*|\geq |X_i^*\cup  X_{i+1}^*\cup T^*|$,  and thus  $|Z|\geq |X_i^*\cup X_{i+1}^*|$. So, the following hold:

\vspace{-0.25cm}
\begin{enumerate}\itemsep=0pt
  \item[$(a)$] For each $i\in [5]$, since $Z\cup A_i\cup A_{i+1}$ is a larger clique than $X_i^*\cup A_i$, we have ${\mathbb{W}}_i=\es$.
  \item[$(b)$] For each $i\in [5]$, since $Z\cup A_{i+2}\cup A_{i-2}$ is a larger clique than $X_i^*\cup X_{i+1}^*\cup A_{i-2}$, we have $X_i^*\cup X_{i+1}^*\cup A_{i-2}\notin {\cal M}$.
 \item[$(c)$] If there is a $T$-clique $T_1$  such that $X_i^*\cup X_{i+1}^*\cup T_1 \in {\cal M}$, then since $|Z|\geq |X_i^*\cup X_{i+1}^*|$, we have $Z\cup T_1 \in {\cal M}$; so $|T_1|\geq 2$ (by \cref{ZXTmax}).
 \end{enumerate}

\vspace{-0.25cm}
\noindent Now,   by \cref{ZXTmax} and $(c)$, $L\cup L'$ meets each maximum clique of $G$ in $G[X\cup Z\cup T]$ twice, and the other maximal cliques in $G[X\cup Z\cup T]$ once.   So we let $S_1:=R_{A_1}\cup R_{A_3}$, $S_2:= R_{A_2}\cup R_{A_5}\cup L$,  and $S_3:=R_{A_4}\cup L'$.  Then, by $(a)$, $(b)$, \cref{XiZ-O3} and \cref{Z-O3},  we see that $S_1,S_2$ and $S_3$ are the required stable sets.

\medskip
\noindent{\bf Case~2} {\it  For each $j\in [5]$, at least one of $X_j, X_{j+2}, X_{j-2}$ is empty. }

\medskip
\no Then there is an index $i\in [5]$ such that $X_i\neq \es$ and $X\sm X_i =\es$ or  $X_{i-1},X_{i}\cup X_{i+1}\neq \es$ and $X\sm (X_{i-1}\cup X_{i}\cup X_{i+1}) =\es$, say $i=1$. By   \cref{ZTX}, \cref{xi-T-Ai+2} and by \ref{TXi-nonadj-last}, we may assume that $A_3$ is a clique.
By \cref{ZTX'}, any $T$-clique $T^*$ that is anticomplete to $X_1$, is complete to an $X_2$-clique or an $X_5$-clique (or to both, if $X_2\cup X_5\neq \es$, by \ref{TXi-nonadj-last}); so  by \cref{xi-T-Ai+2}, $A_4$ is a clique (if $X_2\neq \es$), and $|T'|\geq 2$ (by \cref{ZXTmax}).  See Figure~\ref{fig:Z}:$(c)$.
Now if ${\mathbb{W}}_1\neq \es $, then by \ref{3-XiXi+1-Ai-2}, ${\mathbb{W}}_2\cup {\mathbb{W}}_5=\es$, and we let $S_1:=R_{A_2}\cup R_{A_5}\cup R_{X_1}$, $S_2:=R_{A_1}\cup R_{A_4}\cup L$ and $S_3:=R_{A_3}\cup L'$, and if ${\mathbb{W}}_1=\es$,  then we let $S_1:=R_{A_1}\cup R_{X_2}\cup R_{X_5}$, $S_2:=R_{A_2}\cup R_{A_4}\cup L$ and $S_3:=R_{A_3}\cup R_{A_5}\cup L'$. Then, as earlier by using \cref{XiZ-O3}, it is not hard to verify that
$S_1,S_2$ and $S_3$ are the desired stable sets.

 This completes the proof of Theorem~\ref{thm-5wheel}. \hfill{$\Box$}

\subsection{($P_5$,\,wheel)-free atoms that contains an induced $C_5$}

Since each $k$-wheel, for $k\geq 6$ has an induced $P_5$, by Theorem~\ref{thm-5wheel},   we consider only ($P_5$,\,wheel)-free atoms. Let $G$ be a connected ($P_5$,\,wheel)-free atom that   contains an induced $C_5$,
 say  $v_1$-$v_2$-$v_3$-$v_4$-$v_5$-$v_1$.     Then we define the sets $A$, $X$, $Y$, $Z$ and $T$ as in Section~\ref{sec-C5} with $v_i\in A_i$, for each $i$, and we use the facts \ref{part}--\ref{TXi-nonadj-last} shown in Section~\ref{sec-C5}.   Let $\cal M$ be the set of   maximum cliques of $G$.     Since $G$ is $5$-wheel-free, clearly $Z=\es$.  Thus, if $T\neq \es$, then   $G[T]$ is $P_3$-free (by \ref{5WF-T-P3free}), and recall that by \ref{xi-compT} and \ref{yi-compT}, each vertex in $X\cup Y$ is either anticomplete or good with respect  to $T$. Let $L$ consist of one vertex from each $T$-clique; otherwise let $L:=\es$, and let $L'$ consist of one vertex (which is not in $L$)  from each nontrivial $T$-clique; otherwise let $L':=\es$.
 Moreover, the graph $G$ has some more structural properties, and are given in \cref{prop-xi,lem:XiAi+2maxclqyemp,XnE-Ye-XinacXi+2,lem:C5-Xne-Ye,prop-yi1,F1-YiYi+1Aj-2,prop-yi,clq-ABCD,lem:C5-YneXe} below.

\begin{lemma} \label{prop-xi}
For   $i\in [5]$, the following properties hold:
\begin{lemmalist}
 \item \label{XiXi+1NEclq} Suppose that $X_i$ and $X_{i+1}$ are nonempty. If there is a vertex $p\in A_{i-2}$ that is complete to $X_{i}\cup X_{i+1}$, then $X_i\cup X_{i+1}$ is a clique.
\item \label{XiXi+2-11clq}  Suppose $K$ is an $X_i$-clique and $K'$ is an $X_{i+2}$-clique. Then  $K$ is  complete to $K'$ or $K$ is anticomplete to $K'$.
\item \label{XiAi+2Xi+2Ai}
Let $K$ be an $X_i$-clique, and let $K'$ be an $X_{i+2}$-clique such that $K$ is anticomplete to $K'$. Then either $K$ is complete to $A_{i+2}$ or $K'$ is complete to $A_i$.
\item \label{Xi+1-XiXi+2}
 If $X_{i+1}\neq \es$, then $X_i$ is anticomplete to $X_{i+2}$.
\item \label{WF2-Yemp-TxiXi+1Xi-1}
	No vertex in $T$ has neighbors in three consecutive $X_i$'s.
 \end{lemmalist}
\end{lemma}
 \no{\it Proof of Lemma~\ref{prop-xi}}.~$(i)$:~If $X_i$ and $X_{i+1}$ are cliques, then by \ref{XiXi+1}, the assertion holds. So, up to symmetry, suppose that there are  nonadjacent vertices in $X_i$, say $x$ and $x'$. Let $x''\in X_{i+1}$. Then by \ref{XiXi+1}, $x''$ is complete to $\{x,x'\}$. Also, by our assumption, $p$ is complete to  $\{x, x',x''\}$. Moreover, by \ref{nadj-x-cn},  $x$ and $x'$ have a common neighbor in $A_{i+2}$, say $q$. Now $\{x,q,x',x'',p\}$ induces a $4$-wheel, a contradiction. So $X_i$ is a clique, and by \ref{XiXi+1}, $X_i\cup X_{i+1}$ is a clique. This proves \cref{XiXi+1NEclq}. $\Diamond$

\smallskip
\no{$(ii)$}:~It is enough to show that if a vertex in $K$ has a neighbor in $K'$, then it is complete to $K'$. Suppose not.  Then there are vertices $u\in K$ and $v,w\in K'$ such that $uv,vw\in E$ and $uw\notin E$. If $v$ and $w$ have a common neighbor in $A_i$, say $p$,   then for any neighbor of $u$  in $A_{i+2}$, say $q$,  by \ref{XiAi}, $\{p,u,q,w,v\}$ induces a $4$-wheel, a contradiction. So we may assume that $v$ and $w$ do not share a common neighbor in $A_i$.  So by the definition of $X_{i+2}$ and \ref{xAi+2Ai-2}, both $v$ and $w$ are complete to $A_{i-1}$. Also there is a vertex $r\in A_i$ such that $rv\in E$ and $rw\notin E$.   But then for any neighbor of $u$ in $A_{i+2}$, say  $a$, by \ref{XiAi}, $\{u,r,v_{i-1},w,a,v\}$ induces a $5$-wheel, a contradiction. This proves \cref{XiXi+2-11clq}. $\Diamond$

\smallskip
\no{$(iii)$}:~Suppose there are vertices, say $x\in K$, $a\in A_{i+2}$, $a'\in A_i$ and $x'\in K'$ such that $xa, x'a'\notin E$. Let $a''\in A_{i+2}$ be a neighbor of $x$. Then by \ref{Ai+2-Ai-2-K-xi}, $aa''\notin E$. But then by \ref{XiAi},  $a$-$x'$-$a''$-$x$-$a'$ is a $P_5$, a contradiction. This proves \cref{XiAi+2Xi+2Ai}. $\Diamond$

\smallskip
\no{$(iv)$}:~Let $x\in X_i$ and $x'\in X_{i+2}$, and suppose $x,x'$ are adjacent. Let $u\in X_{i+1}$. By \ref{xAi+2Ai-2}, we may assume that $u$ is complete to $A_{i-2}$.   Now pick a neighbor of $x$ in $A_{i+2}$, say $p$,  and a neighbor of $x$ in $A_{i-2}$, say $q$. Then by \ref{XiAi} and \ref{XiXi+1},  $\{q,u,x',p,x\}$ induces a $4$-wheel, a contradiction. This proves \cref{Xi+1-XiXi+2}. $\Diamond$

\smallskip
\no{$(v)$}:~Suppose there is a  vertex, say $t\in T$ which has neighbors, say $x_1\in X_{1}$, $x_2\in X_2$ and $x_3\in X_{3}$. By \cref{Xi+1-XiXi+2}, $x_1x_3\notin E$. Pick any $a\in A_{4}$ and $a'\in A_{5}$. Then by \ref{xi-pure}, $x_1a,x_2a',x_3a'\in E$, and then $\{t,x_1,a,a',x_3,x_2\}$ induce a $5$-wheel, a contradiction. This proves  \cref{WF2-Yemp-TxiXi+1Xi-1}. \hfill{$\Box$}

\begin{lemma}\label{lem:XiAi+2maxclqyemp}
For $i\in [5]$, let $j,k\in\{i+2,i-2\}$ and $j\neq k$, and let $G[X_i\cup A_{i+2}\cup A_{i-2}]\cong H$. Then the following hold.
\begin{lemmalist}
\item \label{Xi-O3}
  If $M$ is a maximum clique in $H$ such that $M\cap A_{i+2}\neq \es$ and $M\cap A_{i-2}\neq \es$, then $R_{A_{i+2}}\cup R_{A_{i-2}}$ meets $M$ twice.
 \item \label{lemmaXiAi+2maxclq} Let $X^*\subseteq X_i$ be a nonempty clique. If every vertex in $A_{j}$ has a nonneighbor in $X^*$, then $A_{k}$ is a clique.
 \item\label{lemmaXiAi+2meettwice}  If $M$ is a maximum clique in $H$ with  $M\cap A_{j}=\es$, then $M\cap X_i\neq \es$, $A_{k}$ is a clique, and $M\cap A_{k}=A_{k}$. Moreover, $R_{X_i}\cup R_{A_{i+2}}\cup R_{A_{i-2}}$ meets each maximum clique in $H$ twice, and $R_{A_{i+2}}\cup R_{A_{i-2}}$ meets each maximal clique in $H$ at least once.
 \item \label{lemmaXiAi+2maxclqyemp}
	If $Y=\es$, and if there is a maximum clique $M$ in $H$ with $|M|=\omega(G)$ and $M\cap A_{i-2}=\es$ (or $M\cap A_{i+2}=\es$), then $G$ is a nice graph.
\end{lemmalist}
\end{lemma}
\no{\it Proof of \cref{lem:XiAi+2maxclqyemp}}.~$(i)$: If $M\cap X_i=\es$, then, by \ref{Ai-P3free}, clearly the assertion holds; so assume that $M\cap X_i\neq \es$. Let $K$ be an $A_{i+2}$-clique such that $M\cap A_{i+2}\subseteq K$. We claim that $M\cap A_{i+2}=K$. Suppose not, and let $b\in K\sm (M\cap A_{i+2})$. Since $K$ is a clique, $b$ is complete to $M\cap A_{i+2}$. By \ref{Ai+2-Ai-2-K-xi}, $M\cap X_i$ is complete to $b$. By the definition of $A$, $b$ is complete to $M\cap A_{i-2}$. So $b$ is complete to $M$, and hence $M\cup \{b\}$ is a larger clique in $G[X_{i}\cup A_{i+2}\cup A_{i-2}]$, a contradiction; so $M\cap A_{i+2}$ is an $A_{i+2}$-clique. By \ref{Ai-P3free}, $R_{A_{i+2}}$ contains a vertex from each $A_{i+2}$ clique, and $R_{A_{i-2}}$ contains a vertex from each $A_{i-2}$ clique, we see that $R_{A_{i+2}}\cup R_{A_{i-2}}$ meets $M$ twice.  This proves \cref{Xi-O3}. $\Diamond$

\smallskip
\no $(ii)$:~Suppose that $i=1$, $j=4$, and there are nonadjacent vertices in $A_3$, say $a,a'$. Since $v_4\in A_4$,  $v_4$ has a nonneighbor in $X^*,$ say $x$. Let $p$ be a neighbor of $x$ in $A_4$, and let  $x'$ be a nonneighbor of $p$ in $X^*.$  Then, by \ref{xAi+2Ai-2}, $\{x,x'\}$ is complete to $\{a,a'\}$, and then $\{p,a,x',a',x\}$ induces a $4$-wheel, a contradiction. This proves  \cref{lemmaXiAi+2maxclq}. $\Diamond$

\smallskip
\no{$(iii)$}:~To prove the first assertion, we let $j=i-2$. Since $A_{i-2}$ is complete to $A_{i+2}$, clearly $M\cap X_i\neq \es$. Since $M\cap A_{i-2}=\es$, every vertex in  $A_{i-2}$ has a nonneighbor in $M\cap X_i$, and hence, by \cref{lemmaXiAi+2maxclq}, $A_{i+2}$ is a clique. Then, by \ref{Ai+2-Ai-2-K-xi}, $M\cap X_i$ is complete to $A_{i+2}$; so $M\cap A_{i+2}=A_{i+2}$. To prove the second assertion, let $M'$ be a maximum clique in $H$. By \cref{Xi-O3}, we may assume that one of $M'\cap A_{i+2}= \es$, $M'\cap A_{i-2}=\es$. If $M'\cap A_{i-2}=\es$, then by the first assertion, since $A_{i+2}$ is a clique,  $M'=X^*\cup A_{i+2}$, where $X^*$ is an $X_i$-clique. Thus $R_{X_i}\cup R_{A_{i+2}}$ meets $M'$ twice, and $R_{A_{i+2}}$ meets $M'$ at least once. Likewise, if $M'\cap A_{i+2}=\es$, then $R_{X_i}\cup R_{A_{i-2}}$ meets $M'$ twice, and $R_{A_{i-2}}$ meets $M'$ at least once. This proves  \cref{lemmaXiAi+2meettwice}.  $\Diamond$

\smallskip
\no{$(iv)$}:~To prove the assertion, we let $i=1$, and suppose that $M\cap A_4=\es$.  As shown in the proof of second assertion of \cref{lemmaXiAi+2meettwice},  $M=X^*\cup A_3$, where $X^*$ is an $X_1$-clique. Let $x\in X^*$, and let $a\in A_4$ be a neighbor of $x$. Then $a$ has a non-neighbor in $X^*$, say $x'$. Then:

\vspace{-0.25cm}
\begin{enumerate}[label=(\alph*)]\itemsep=0pt
\item For any $p\in X_5$, by \ref{XiXi+1} and \ref{Ai+2-Ai-2-K-xi},    $M \cup \{p\}$ is a clique, a contradiction; so $X_5=\es$.
\item If there is a vertex $p\in X_3$,  for any neighbor of $p$ in $A_5$, say $q$, since $p$-$q$-$a$-$x$-$x'$ is not a $P_5$, $p$ is adjacent to one of $x,x'$, then, by \cref{XiXi+2-11clq}, $p$ is complete to $X^*$, and then, by \ref{XiAi}, $M\cup \{p\}$ is a clique, a contradiction; so $X_3=\es$.
\item Suppose there is a vertex  $p\in X_4$. Then for any neighbor of $p$ in $A_2$, say $q$, $q$-$p$-$a$-$x$-$x'$ is not a $P_5$,  $p$ is adjacent to one of $x,x'$. Let $K$ be the $X_4$-clique containing $p$. Then, by \cref{XiXi+2-11clq}, $K$ is complete to $X^*$, and then, by \ref{KcompK'}, $X^*$ is complete to exactly one $A_4$-clique, say $K'$. Then since $M\cup K'$ is a clique, a contradiction. So $X_4=\es$.
\item If there are adjacent vertices, say $t\in T$ and $x_2\in X_2$, and if $K$ is the $X_2$-clique containing $x_2$, and $Q$ is the $A_4$-clique containing $a$, then by \ref{xi-pure}, $x_2$ is complete to $A_4$, and then since $N(K)\cap Q\neq \es$ and $N(X^*)\cap Q\neq \es$, by \ref{G7-clm1}, $ax'\in E$, a contradiction; so $X_2$ is anticomplete to $T$.
\item If $M'$ is a maximal clique in $G$ such that $M'\cap T\neq \es$, then since $G$ is an atom, by $(d)$, $M'\cap X_1\neq \es$, then, by \ref{xi-pure}, for any $A_{3}$-clique  $D$, and any $A_4$-clique $D'$, $|(M'\cap X_1)\cup D\cup D'|\leq M$. Hence $|M'\cap T|\geq 2$.
\end{enumerate}
Now by \ref{XiAi},  \ref{xi-compT} and \cref{lemmaXiAi+2meettwice}, the sets $S_1:=R_{A_2}\cup R_{A_5}\cup R_{X_1}$, $S_2:=R_{A_1}\cup R_{A_3}\cup R_{X_2}\cup L$, and $S_3:=R_{A_4}\cup L'$ are the required stable sets. So $G$ is nice. This proves \cref{lemmaXiAi+2maxclqyemp}. \hfill{$\Box$}

\begin{lemma}\label{XnE-Ye-XinacXi+2}
If $Y$ is empty, and if there is an $i\in [5]$ such that $X_i$ is not anticomplete to $X_{i+2}$, then $G$ is nice.
\end{lemma}
\no{\it Proof of \cref{XnE-Ye-XinacXi+2}}. We may assume that $i=1$. Then there are vertices $x_1\in X_1$ and $x_3\in X_3$ such that $x_1x_3\in E$. Then by \cref{Xi+1-XiXi+2}, $X_2=\es$; so $X=X_1\cup X_3\cup X_4\cup X_5$. Let $Q_1$ be the $X_1$-clique containing $x_1$, and let $Q_3$ be  the $X_3$-clique containing $x_3$.
Then by \cref{XiXi+2-11clq} and  \ref{KcompK'},   $Q_1$ is complete to $Q_3$, and anticomplete to $(X_3\sm Q_3)\cup X_4$, $Q_3$ is anticomplete to $(X_1\sm Q_1)\cup X_5$, and  $X_1\sm Q_1$ is anticomplete to $X_3\sm Q_3$.  By \ref{KcompK'}, let $A_1^*$ be the $A_1$-clique such that $Q_3$ is complete to $A_1^*$, and anticomplete to $A_1\sm A_1^*$, and  let $A_3^*$ be the $A_3$-clique  such that $Q_1$ is complete to $A_3^*$, and anticomplete to $A_3\sm A_3^*$. By \ref{nadj-x-cn} and \ref{x-good}, $X_1$ is complete to $A_3^*$, and $X_3$ is complete to $A_1^*$.

Note that any maximal clique containing at least one vertex from each $X_1$ and $X_3$ is either $A_1^*\cup Q_1\cup Q_3$ or $A_3^*\cup Q_1\cup Q_3$. By \ref{KcompK'}, any maximal clique containing at least one vertex from each $X_1$ and $X_4$ is of the form $D_1\cup X_1^*\cup X_4^*$ or $D_4\cup X_1^*\cup X_4^*$, where $D_1$, $D_4$, $X_1^*$ and $X_4^*$ are $A_1$, $A_4$, $X_1$ and $X_4$-cliques respectively, and $X_1^*\neq Q_1$. Also, any maximal clique containing at least one vertex from each $X_3$ and $X_5$ is of the form $D_3\cup X_3^*\cup X_5^*$ or $D_5\cup X_3^*\cup X_5^*$, where $D_3$, $D_5$, $X_3^*$, $X_5^*$ are $A_3$, $A_5$, $X_3$ and $X_5$-cliques respectively, and $X_3^*\neq Q_3$.

By \cref{lemmaXiAi+2maxclqyemp}, we may assume that each maximal clique of $G$ in $G[X_i\cup A_{i+2}\cup A_{i+2}]$  has nonempty intersection with both $A_{i+2}$ and $A_{i-2}$; and by \cref{lemmaXiAi+2meettwice}, $R_{A_{i+2}}\cup R_{A_{i-2}}$ meets rest of the maximal cliques in $G[X_i\cup A_{i+2}\cup A_{i+2}]$ at least once.

First  suppose that $T=\es$.    Also assume that  $Q_1$  is either complete or anticomplete to every $A_4$-clique,  and  $Q_3$ is either complete or anticomplete to  every $A_5$-clique.
 Now suppose there is an $A_2$-clique, say $D_2$, such that either $A_1^*\cup D_2\in {\cal M} $ or $A_3^*\cup D_2\in {\cal M}$. Up to relabelling, we may assume that $A_1^*\cup D_2\in {\cal M}$.  Then since $|A_1^*\cup D_2|\geq |A_1^*\cup Q_1\cup Q_3|$, we have $|D_2|>|Q_1|$.  Further, we have the following:

\vspace{-0.25cm}
\begin{enumerate}\itemsep=0pt

\item[$(a)$] Any maximal clique that contain at least one vertex from each $A_4$ and $Q_1$ is of the form $A_3^*\cup Q_1\cup D_4$, where $D_4$ is an $A_4$-clique.

\item[$(b)$] For any $A_1$-clique $D_1$, since $|D_2|>|Q_1|$,  we have $D_1\cup Q_1\notin {\cal M}$.

\item[$(c)$] If $X_5\neq \es$,  since $X_5$ is anticomplete to $Q_3$, by \ref{XiAi-2clqcor}, each $X_5$-clique is either complete or anticomplete to an $A_3$-clique. So for any $X_5$-clique  $X_5^*$ which is anticomplete to $A_3^*$, by \ref{xAi+2Ai-2}, $X_5^{*}$ is complete to $D_2$, and   $|Q_1\cup X_5^{*}|< |D_2\cup X_5^{*}|$  which implies that $Q_1\cup X_5^{*}\notin {\cal M}$. Moreover, for any  $X_5$-clique $X_5^{**}$ which is complete to $A_3^*$, any maximal clique that contain at least one vertex from each $X_5^{**}$ and $Q_1$ is of the form $A_3^*\cup Q_1\cup X_5^{**}$.
\end{enumerate}

\vspace{-0.25cm}
\noindent By $(a)$, $(b)$ and $(c)$, it is not hard to verify that $S_1:=R_{A_2}\cup R_{X_1\sm Q_1} \cup R_{X_3}$, $S_2:=R_{A_3}\cup R_{A_5}\cup R_{X_4}$, and $S_3:=R_{A_1}\cup R_{A_4}\cup R_{X_5}$ are the required stable sets.  So we assume that for any $A_2$-clique  $D_2$, $A_1^*\cup D_2, A_3^*\cup D_2\notin {\cal M}$. Next we claim the following:
\begin{claim}\label{W1A1W3A3}
Either for each $W\in {\mathbb{W}}_1$, $W\cap A_1^*=\es$ or for each $W'\in {\mathbb{W}}_3$, $W'\cap A_3^*=\es$.
\end{claim}
\no{\it Proof of \cref{W1A1W3A3}}.~Suppose   there is an $X_1$-clique $K$ such that $K\cup A_1^* \in {\mathbb{W}}_1$, and  there is an $X_3$-clique $K'$ such that $K'\cup A_3^* \in {\mathbb{W}}_3$.  Note that $K\neq Q_1$ and $K'\neq Q_3$. Then $K$ is anticomplete to $K'$. Let $D_5$ be an $A_5$-clique such that  $N(K')\cap D_5\neq \es$. Then, by \ref{XiAi-2clqcor}, $K'$ is complete to $D_5$. Now $|A_1^*\cup K|\geq |A_1^*\cup D_5\cup K'|$, and so $|K|>|K'|$. Then $A_3^*\cup K$ is a clique, and  $|A_3^*\cup K|> |A_3^*\cup K'|$ which is a contradiction. This proves \cref{W1A1W3A3}. $\Diamond$

\smallskip
By \cref{W1A1W3A3}, we assume, up to symmetry, that for each  $W\in {\mathbb{W}}_1$, we have $W\cap A_1^*=\es$.
Now if for each $A_5$-clique $D_5$, $D_5\cup A_1^*\notin {\cal M}$, then clearly $S_1:=R_{A_1\sm A_1^*}\cup R_{A_4}\cup R_{Q_3}\cup R_{X_5}$, $S_2:=  R_{A_2}\cup R_{X_1}\cup R_{X_3\sm Q_3}$ and $S_3:=R_{A_3}\cup R_{A_5}\cup R_{X_4}$ are the required stable sets. So suppose that there is an $A_5$-clique $D_5$ such that $D_5\cup A_1^*\in{\cal M}$. Then since $Q_3\cup D_5\cup A_1^*$ is not a clique,  $Q_3$ is anticomplete to $D_5$. Then, by \ref{xAi+2Ai-2}, $Q_3$ is complete to $A_1$;  so $A_1=A_1^*$, and hence ${\mathbb{W}}_1=\es$. Also, if $X_5\neq \es$, then since $Q_3$ is anticomplete to $X_5$, by \cref{XiAi+2Xi+2Ai}, $X_5$ is complete to $A_3$. Thus, by \ref{G7-clm1}, any maximum clique containing at least one vertex from each $X_1$ and $X_5$ is of the form $D_3\cup X_1^*\cup X_5^*$, where $X_1^*$, $X_5^*$ and $D_3$ are $X_1$, $X_5$ and $A_3$-clique respectively. Now we let $S_1:=R_{A_2}\cup R_{A_4}\cup R_{X_3}$, $S_2:=R_{A_3}\cup R_{A_5}\cup R_{X_4}$ and $S_3:=R_{A_1}\cup R_{X_5}$, and we conclude that $S_1$, $S_2$ and $S_3$ are the required stable sets.
So suppose that, up to relabelling, there is an $A_5$-clique, say $D_5$, such that $Q_3$ is neither complete nor anticomplete to $D_5$. Then since $(X_1\sm Q_1)\cup X_5$ is anticomplete to $Q_3$, by \ref{XiAi-2clqcor}, $(X_1\sm Q_1)\cup X_5=\es$. So $X_1=Q_1$ and $X_1$ is anticomplete to $X_4$ (by \ref{KcompK'}). Now   we let $S_1:=R_{A_5}\cup R_{X_1}\cup R_{X_4}$, $S_2:=R_{A_2}\cup R_{A_4}\cup R_{X_3}$ and $S_3:=R_{A_1}\cup R_{A_3}$. Then clearly $S_1, S_2$ and $S_3$ are the required stable sets.

\begin{figure}[h]
\centering
 \includegraphics[height=4cm, width=14cm]{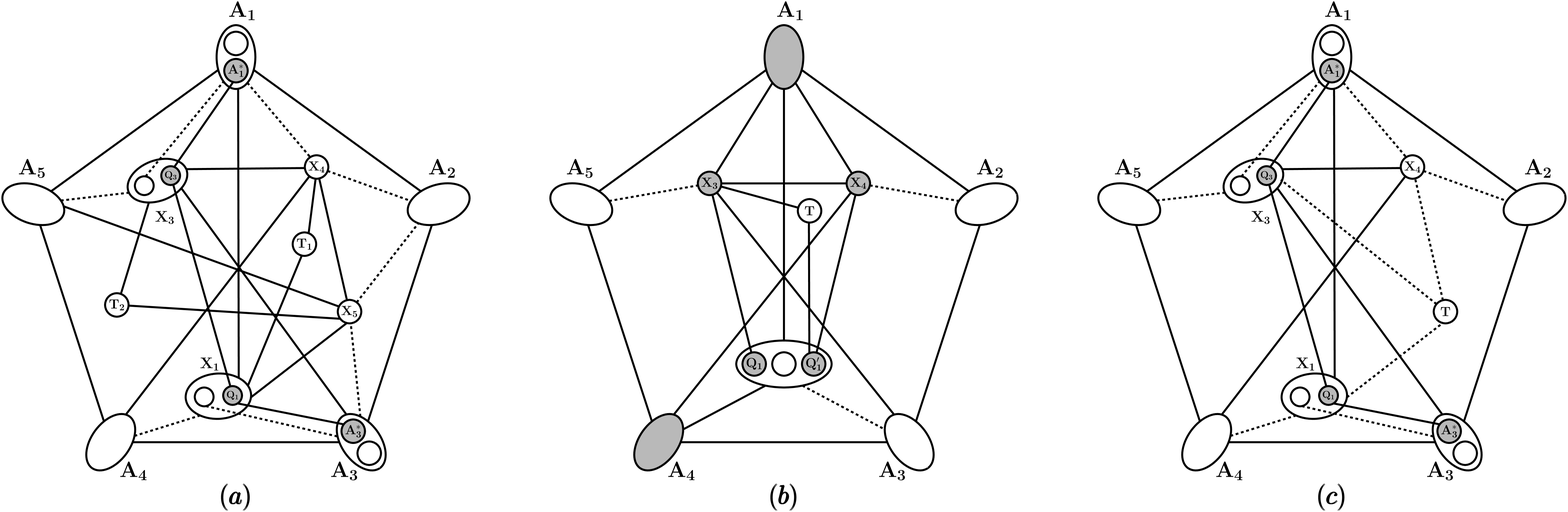}
\caption{Sketch of the graph $G$ in   \cref{XnE-Ye-XinacXi+2} when $T\neq \es$, and: $(a)$  $X_4,X_5\neq \es$. $(b)$  $X_1$ is not anticomplete to $X_4$. $(c)$ $X_1$ is  anticomplete to $X_4$.}\label{fig:xT}
\end{figure}

So, we may assume that $T\neq \es$.  (We refer to Figure~\ref{fig:xT} for a sketch of the graph $G$.) Recall that each vertex in $X$ is either anticomplete or good with respect to $T$.  First suppose $X_4,X_5\neq \es$. So, by \cref{Xi+1-XiXi+2}, $X_1$ is anticomplete to $X_4$, and $X_3$ is anticomplete to $X_5$. Let $T_1$ denote the union of $T$-cliques which are complete to $X_1\cup X_4$, and anticomplete to $X_3\cup X_5$, and let $T_2$ denote the union of $T$-cliques which are complete to $X_3\cup  X_5$, and anticomplete to $X_1\cup X_4$. Clearly $T_1\cap T_2=\es$. Moreover:

\begin{claim}\label{lc-Tpart}
 $T=T_1\cup T_2$.
\end{claim}
\no{\it Proof of \cref{lc-Tpart}}.~Let $t\in T$, and let $T'$ be the $T$-clique containing $t$. Since every vertex in $T$ has a neighbor in $X$,   first assume that $t$ has a neighbor in $X_1\cup X_4$. Since $X_1$ is anticomplete to $X_4$, by \ref{TXi-nonadj-last}, $t$ is complete to $X_1\cup X_4$. So by \cref{WF2-Yemp-TxiXi+1Xi-1}, $t$ is anticomplete to $X_5$, and since $X_3$ is anticomplete to $X_5$, by \ref{TXi-nonadj-last}, $t$ is anticomplete to $X_3$.  Thus, by \ref{xi-compT}, $T'$ is complete to $X_1\cup X_4$, and anticomplete to $X_3\cup X_5$, and so $T'\in T_1$. Similarly, if $t$ has a neighbor in $X_3\cup X_5$, then $T'$ is complete to $X_3\cup X_5$, and anticomplete to $X_1\cup X_4$, and so $T'\in T_2$. This proves \cref{lc-Tpart}. $\Diamond$

\medskip
  Since $T\neq \es$, by \ref{xi-pure}, for $j\in\{3,4,5\}$,  either $X_j$ is complete to $A_{j-2}$ or $X_{j+1}$ is complete to $A_{j-2}$; so any maximal clique containing at least one vertex from each $X_j$ and $X_{j+1}$ must be complete to  $A_{j-2}$ (by \ref{G7-clm1}). Since at least one of $X_4$ and $X_5$ is complete to $A_2$, by \ref{G7-clm1} and \ref{3-XiXi+1-Ai-2},  one of ${\mathbb{W}}_4$ and ${\mathbb{W}}_5$ is empty. Moreover, if $T^*\cup X^*\in {\cal M}$ for a $T$-clique $T^*$, and an $X_i$-clique $X^*$, where $i\in \{1,3,4,5\}$, then by \ref{xi-pure}, for any $p\in A_{i+2}$ and $q\in A_{i-2}$, $|T^*\cup X^*|\geq |X^*\cup \{p,q\}|$, and thus $|T^*|\geq 2$.   Now, if ${\mathbb{W}}_5=\es$, then we let $S_1:=R_{A_5}\cup R_{X_1}\cup R_{X_4}\cup (L\cap T_2)$, $S_2:=R_{A_2}\cup R_{A_4}\cup R_{X_3}\cap (L\cap T_1)$ and $S_3:=R_{A_1}\cup R_{A_3}\cup  (L'\cap T_2)$, and if  ${\mathbb{W}}_4=\es$, then we let  $S_1:=R_{A_4}\cup R_{X_3}\cup R_{X_5}\cup (L\cap T_1)$, $S_2:=R_{A_2}\cup R_{A_5}\cup R_{X_1}\cup (L\cap T_2)$, $S_3:=R_{A_1}\cup R_{A_3}\cup (L'\cap T_1)$. Then we observe that $S_1, S_2$, and $S_3$ are the required stable sets.

Next we assume that one of $X_4$ and $X_5$ is empty, say $X_5=\es$. First suppose that $X_1$ is not anticomplete to $X_4$. So there are vertices $x_1'\in X_1$,  $x_4\in X_4$  such that $x_1'x_4\in E$. So $X=X_1\cup X_3\cup X_4$. Let $Q_1'$ be the $X_1$-clique containing $x_1'$, and let $Q_4$ be  the $X_4$-clique containing $x_4$.
Then by \ref{KcompK'} and \cref{XiXi+2-11clq},  $Q_1\neq Q_1'$,  $Q_1'$ is complete to $Q_4$, $Q_1'$ is anticomplete to $X_4\sm Q_4$, $Q_4$ is anticomplete to $X_1\sm Q_1'$,  $X_1\sm Q_1'$ is anticomplete to $X_4\sm Q_4$.
By \ref{KcompK'},  let $A_1^{**}$ be  the $A_1$-clique such that $Q_4$ is complete to $A_1^{**}$ and anticomplete to $A_1\sm A_1^{**}$, and  let $A_4^*$ be the $A_4$-clique  such that $Q_1'$ is complete to $A_4^*$, and anticomplete to $A_4\sm A_4^*$.  By \ref{nadj-x-cn} and \ref{x-good}, $X_4$ is complete to $A_1^{**}$. By \ref{TXi-nonadj-last}, each vertex in $T$ has a neighbor in $X_1$. Further we claim the following:

\begin{claim}\label{Tjust}
   $A_1$,  $X_3$ and $X_4$ are cliques.  Moreover, $T$ is complete to  exactly one of $X_3$ and $X_4$.
\end{claim}
\no{\it Proof of \cref{Tjust}}.~We first show that, if $x_3$ is not anticomplete to $T$, then $x_4$ is anticomplete to $T$, and vice versa. Suppose not, and let $t,t'\in T$ be such that $x_3t,x_4t'\in E$. If $x_4t\in E$, then, by \ref{TXi-nonadj-last}, $x_1t\in E$, and then by \ref{xi-pure}, for any $a\in A_1$, $\{a, x_1,t,x_4,x_3\}$ induces a $4$-wheel, a contradiction; so $x_4t\notin E$. Likewise, $x_3t'\notin E$. Also, by \ref{xi-compT}, $tt'\notin E$, and by \ref{TXi-nonadj-last}, $x_1t',x_1't\in E$ and $x_1t,x_1't'\notin E$. But then $t$-$x_1'$-$x_4$-$t'$-$x_1$ is a $P_5$, a contradiction. By symmetry, we may assume that $x_3$ is not anticomplete to $T$. Then $x_4$ is anticomplete to $T$. Then, by \ref{TXi-nonadj-last}, $T$ is anticomplete to $X_1\sm Q_1'$. Since each vertex in $T$ has a neighbor in $X_1$, each vertex in $T$ has a neighbor in $Q_1'$. So by \ref{xi-compT} and \ref{TXi-nonadj-last}, $T$ is complete to $Q_1'\cup Q_3$. By \ref{xi-pure}, $A_1=A_1^*=A_1^{**}$ is a clique. So by \cref{XiXi+1NEclq}, $X_3\cup X_4$ is a clique, and hence   $X_3=Q_3$ and $X_4=Q_4$ are cliques.  Since     $T$ is anticomplete to $X_1\sm Q_1'$, by \ref{TXi-nonadj-last}, $T$ is anticomplete to $Q_4=X_4$. This proves \cref{Tjust}. $\Diamond$

\smallskip
By \cref{Tjust},    we may assume that $T$ is complete to $X_3(=Q_3)$, and anticomplete to $X_4(=Q_4)$.
Then by \ref{TXi-nonadj-last}, it follows that, $T$ is complete to $Q_1'$ (and anticomplete to $X_1\sm Q_1'$).   Then by \ref{xi-pure}, $Q_1'$ is complete to $A_4$, and hence $A_4=A_4^*$ is a clique. Since $Q_1'\cup Q_4\cup A_4$ is a larger clique than $Q_4\cup A_4$, we conclude that ${\mathbb{W}}_4=\es$, and $R_T\cup R_{X_3}\cup R_{Q_1'}$ meets each maximal clique of $G$ in $G[X\cup T]$ twice. Now we  see that $S_1:=R_{A_2}\cup R_{A_5}\cup R_{X_1}$, $S_2:=R_{A_4}\cup R_{X_3}$, and $S_3:=R_{A_1}\cup R_{A_3}\cup R_T$ are the required stable sets.

Finally we assume that either $X_1$ is anticomplete to $X_4$ or $X_4=\es$. We claim the following:
\begin{claim}\label{Q1Q3just}
Each $T$-clique is  complete to either $Q_1$ or $Q_3$.
\end{claim}
\no{\it Proof of \cref{Q1Q3just}}.~Suppose not. Then there is a $T$-clique, say $T^*$, and vertices $x\in Q_1$, $x'\in Q_3$, and $t,t'\in T^*$ such that $xt,xt'\notin E$. Then by \ref{xi-compT}, $T^*$ is anticomplete to $\{x,x'\}$. So by \ref{KcompK'} and \ref{TXi-nonadj-last}, $T^*$ is anticomplete to $(X_1\sm Q_1)\cup (X_3\sm Q_3)\cup X_4$. Since each vertex of $T$ has a neighbor in $X$, $N(T^*)\cap X\subseteq Q_1\cup Q_3$ which is a clique cutset, a contradiction. This proves \cref{Q1Q3just}. $\Diamond$

\smallskip
Moreover, if there is a $T$-clique $T^*$ such that $N(T^*)\cap (X_4\cup (X_3\sm Q_3))\neq \es$, then by \ref{TXi-nonadj-last}, $T^*$ is complete to $Q_1$, and hence  $T^*$ is complete to $X_4\cup (X_3\sm Q_3)$.   Now we let  $S_1:=R_{A_5}\cup R_{X_1}\cup R_{X_4}$, $S_2:=R_{A_2}\cup R_{A_4}\cup R_{X_3}$, and $S_3:=R_{A_1}\cup R_{A_3}\cup R_{T}$. Then by \cref{Q1Q3just},  we conclude that $S_1, S_2$ and $S_3$ are the desired stable sets. This completes the proof of \cref{XnE-Ye-XinacXi+2}. \hfill{$\Box$}

\begin{figure}[h]
\centering
 \includegraphics[height=3.5cm, width=12cm]{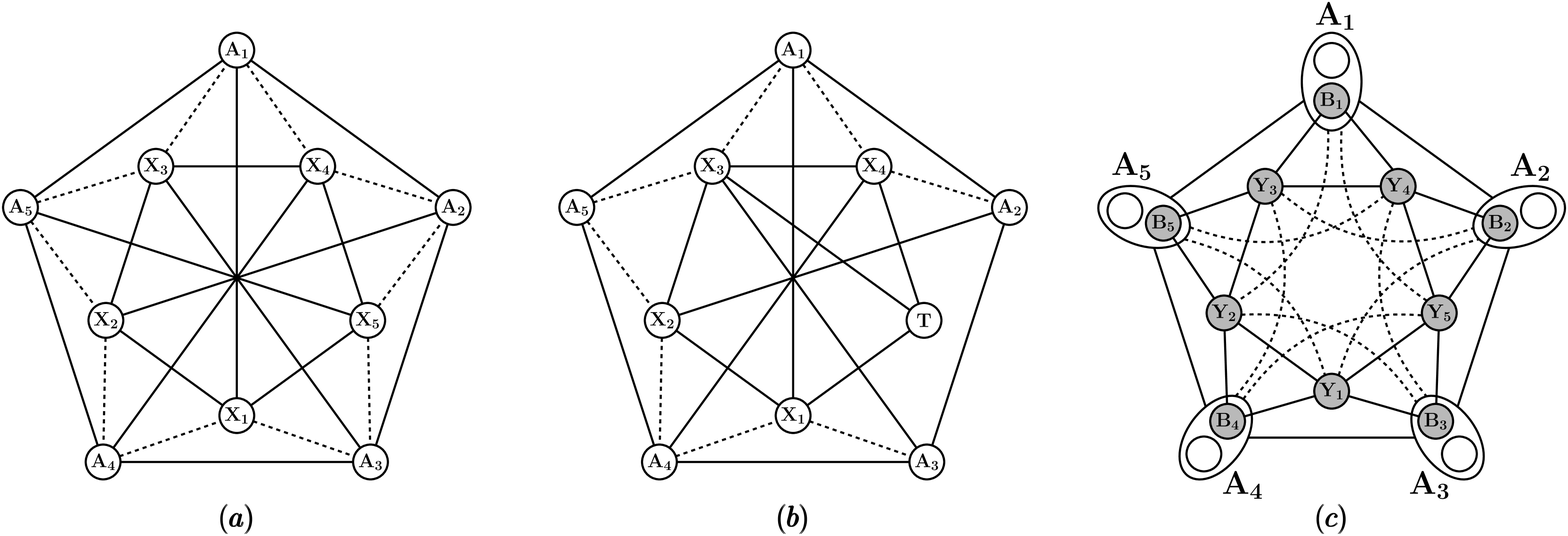}
\caption{Sketch of the graph $G$ in: $(a)$ Lemma~\ref{lem:C5-Xne-Ye} when $T=\es$. $(b)$ Lemma~\ref{lem:C5-Xne-Ye} when $T\neq \es$ and $\ell=5$. $(c)$ Lemma~\ref{lem:C5-YneXe}  when $Y_i$ is anticomplete to $Y_{i+2}\cup Y_{i-2}$, for each $i\in [5]$.}\label{fig:xory}
\end{figure}

\begin{lemma}\label{lem:C5-Xne-Ye}
 If $X$ is nonempty and $Y$ is empty, then $G$ is a nice graph.
\end{lemma}
\no{\it Proof of Lemma~\ref{lem:C5-Xne-Ye}}.~By \cref{lemmaXiAi+2maxclqyemp}, we may assume that each maximal clique of $G$ in $G[X_i\cup A_{i+2}\cup A_{i+2}]$  has nonempty intersection with both $A_{i+2}$ and $A_{i-2}$; and by \cref{lemmaXiAi+2meettwice}, $R_{A_{i+2}}\cup R_{A_{i-2}}$ meets rest of the maximal cliques in $G[X_i\cup A_{i+2}\cup A_{i+2}]$ at least once.
Recall that each vertex in $X$ is either anticomplete or good with respect to $T$.  First suppose that  $X\sm X_1=\es$.  If $T^*$ is a $T$-clique such that $T^*\subset M \in \cal{M}$, then since each vertex in $X_1$ is either anticomplete or good with respect to $T$, $M = T^*\cup X_1^*$ where $X_1^*$ is a subset of some $X_1$-clique. Since $N(T^*)\cap X_1$ is  not a clique cutset, there are nonadjacent vertices in $N(T^*)\cap X_1$. Then by \ref{2nax-Ai+2Ai-2clq} and \ref{xi-pure}, $A_3\cup A_4$ is clique, and so $A_3\cup A_4\cup X_1^*$ is a clique. Hence $|T^*\cup X_1^*|\geq |A_3\cup A_4\cup X_1^*|$, and thus $|T^*|\geq 2$. Then clearly  $S_1:=R_{X_1}\cup R_{A_2}\cup R_{A_5}$, $S_2:=R_{A_1}\cup R_{A_3}\cup L$, and $S_3:=R_{A_4}\cup L'$ are the desired stable sets. Let $J$ denote the set $\{i\in [5]\mid X_i\neq \es\}$, and we may assume that $|J|\ge 2$.
 By \cref{XnE-Ye-XinacXi+2}, we may assume that for each $i\in [5]$, $X_i$ is anticomplete to $X_{i+2}$.
See Figure~\ref{fig:xory}:$(a)$ and Figure~\ref{fig:xory}:$(b)$.   First we claim the following.

\begin{claim}\label{Wiempupd}
There is an index $\ell\in [5]$ such that ${\mathbb{W}}_{\ell}=\es$, and for $p\in \{\ell+1, \ell-1\}$,  $ R_{X_p} \cup R_{A_1}\cup \cdots \cup R_{A_5}$ meets each maximum clique of $G$ in $G[A\cup X_{\ell}\cup X_{p}]$  at least twice.
\end{claim}
\no{\it Proof of \cref{Wiempupd}}.~If there is an index $i\in [5]$ such that $X_i=\es$,  we choose $\ell=i$, and by \cref{Xi-O3}, we are done; so for each $i\in [5]$, $X_i\neq \es$. First suppose there is an index $i\in [5]$, $X_i$ and $X_{i+1}$ is not complete to $A_{i-2}$, say $i=1$. Then by \cref{XiAi+2Xi+2Ai}, $X_4$ is complete to $A_1\cup A_2$, and, up to relabelling,  $X_3$ is complete to $A_5$. Since $X_4$ is complete to $A_1$, by \ref{3-XiXi+1-Ai-2}, one of ${\mathbb{W}}_3$ and ${\mathbb{W}}_4$ is empty.
Now if  ${\mathbb{W}}_3=\es$, then we choose $\ell=3$, otherwise we choose $\ell=4$.  Note that, by \ref{G7-clm1},  for $k\in \{2,3,4\}$, any maximal clique containing at least one vertex from each $X_k$ and $X_{k+1}$ must contain a vertex from $R_{A_{k-2}}$. Thus, by \cref{Xi-O3}, we conclude the proof. So we may assume that for each $i\in [5]$, one of $X_i$ or $X_{i+1}$ is complete to $A_{i-2}$. Then by \ref{3-XiXi+1-Ai-2} and \ref{G7-clm1}, there is an index $k\in [5]$ such that ${\mathbb{W}}_k=\es$.  Since  by \ref{G7-clm1}, for $i\in [5]$, any maximal clique containing at least one vertex from each $X_i$ and $X_{i+1}$ must contain a vertex from $R_{A_{i-2}}$,  by \cref{Xi-O3}, we conclude that  $\ell= k$ is our desired index.  This proves \cref{Wiempupd}. $\Diamond$

\begin{claim}\label{3.1-TX}
	If $T\neq \es$, then there is an index $j\in [5]$ such that  $X_{j-1}=\es$ (so ${\mathbb{W}}_{j-1}=\es$), and each $T$-clique  is complete to  $X_j\cup X_{j+2}\cup X_{j-2}$, and anticomplete to $X_{j+1}$.
\end{claim}
\no{\it Proof of \cref{3.1-TX}}.~Let $t\in T$. Let $T^*$ be the $T$-clique containing $t$. By \cref{F1-TXY2}, there is an $i\in [5]$ such that $N(T^*)\cap X_i\neq \es$. Suppose $N(T^*)\cap X_{i+1} \neq \es$.  Then by \ref{xi-compT}, \ref{TXi-nonadj-last} and \cref{WF2-Yemp-TxiXi+1Xi-1}, we may assume that $X_{i+2}$ and $X_{i-1}$ are empty.  If $X_{i-2}= \es$, then by \ref{xi-pure} and  \cref{XiXi+1NEclq},  $X_i\cup X_{i+1}$ is a clique, and so $N(T^*)\cap (X_i\cup X_{i+1})$ is a clique cutset, a contradiction; so $X_{i-2}\neq \es$. Then by \ref{TXi-nonadj-last}, $T^*$ is complete to $X_i\cup X_{i+1} \cup X_{i-2}$. So we take $j=i-2$ and we are done.
Thus, by \ref{xi-compT}, we may assume that  $N(T^*)\cap (X_{i-1}\cup X_{i+1}) =\es$. We claim that $X_{i-2}\cup X_{i+2}\neq \es$. Suppose not. Since $|J|\geq 2$, we may assume that $X_{i+1}\neq \es$. Then since  $N(T^*)\cap X_i$ is not a clique cutset, there are nonadjacent vertices, say $u,v\in X_i$ such that $u,v\in N(T^*)$. Then by \ref{2nax-Ai+2Ai-2clq} and \ref{xi-pure}, $A_{i-2}$ is a clique, and so $X_i\cup X_{i+1}$ is complete to $A_{i-2}$ (by \ref{Ai+2-Ai-2-K-xi}). But then, by \cref{XiXi+1NEclq}, $X_i$ is a clique, a contradiction. Thus  $X_{i-2}\cup X_{i+2}\neq \es$. Then by \ref{TXi-nonadj-last}, $T^*$ is complete to $X_i\cup X_{i+2} \cup X_{i-2}$. Also, by \cref{WF2-Yemp-TxiXi+1Xi-1}, $T^*$ is anticomplete to $X_{i+1}$. So we take $j=i$. This proves \cref{3.1-TX}. $\Diamond$

\smallskip
By \cref{Wiempupd},  let $\ell\in [5]$ be the index such that ${\mathbb{W}}_{\ell}=\es$, and for $p\in \{\ell+1, \ell-1\}$,  $R_{X_p} \cup R_{A_1}\cup \cdots \cup R_{A_5}$ meets each maximum clique of $G$ in $G[A\cup X_{\ell}\cup X_{p}]$  at least twice. If $T\neq\es$, then we choose $\ell=j-1$ (by \cref{3.1-TX}).
Now we let $S_1:=R_{X_{\ell-1}}\cup R_{A_\ell}\cup R_{X_{\ell+1}}$, $S_2:=R_{A_{\ell+1}}\cup R_{X_{\ell+2}}\cup R_{A_{\ell-2}}\cup R_T$, and $S_3:=R_{A_{\ell+2}}\cup R_{X_{\ell-2}}\cup R_{A_{\ell-1}}$. Clearly  $S_1,S_2$ and $S_3$ are stable sets. By \cref{3.1-TX}, $R_{X_{j-2}}\cup R_{X_{j}}\cup R_{X_{j+1}}\cup R_{X_{j+2}}\cup R_T$ meets each maximum clique of $G$ in $G[X\cup T]$ at least twice. Also, using \cref{Xi-O3} and \ref{G7-clm1}, we see that $ (\cup_{k=1}^5 R_{A_k}) \cup R_{X_{\ell-2}}\cup R_{X_{\ell-1}}\cup R_{X_{\ell+1}}\cup R_{X_{\ell+2}}$ meets each maximum clique of $G$ in  $G[A\cup (X\sm X_{\ell})]$ twice.
So, by \cref{Wiempupd}, we conclude that $S_1\cup S_2\cup S_3$ meets each maximum clique of $G$ at least twice, and other maximal cliques at least once. So $G$ is nice. This completes the proof of \cref{lem:C5-Xne-Ye}. \hfill{$\Box$}

\begin{lemma} \label{prop-yi1}
For   $i\in [5]$, the following properties hold:
\begin{lemmalist}
\item \label{Ai+2-Ai-2-K-yi}
 If $K$ is an $A_{i+2}$-clique (or an $A_{i-2}$-clique), then any vertex in $Y_i$ which has a neighbor in $K$ is complete to $K$.
\item\label{Ai+2-Ai-2-yi}
For  $j\in \{i-2,i+2\}$, each vertex in  $Y_i$ is complete to exactly one $A_j$-clique. 
\item\label{YiYi+2Yi-2comYi+1ancom}
If a vertex in $Y_i$ is not complete to $A_{i-1}$ (or $A_{i+1}$), then it is complete to $A_{i+2}\cup A_{i-2}$, and so $A_{i+2}\cup A_{i-2}$ is a clique.
\item\label{Yi-clq}
 $Y_i$ is a clique.
 \item\label{YiYi+1-clq}
 $Y_i$ is complete to $Y_{i+1}\cup Y_{i-1}$.
\item\label{F1-TXY2} Every vertex in $T$ has a neighbor in $X$.
 \end{lemmalist}
\end{lemma}
\no{\it Proof of \cref{prop-yi1}}.~$(i)$:~The proof of \cref{Ai+2-Ai-2-K-yi} is similar to that of \ref{Ai+2-Ai-2-K-xi}, and we omit the details. $\Diamond$

\smallskip
\no{$(ii)$}:~We may assume, up to symmetry,  that $j=i+2$. Let $y\in Y_i$. By \cref{Ai+2-Ai-2-K-yi}, it is enough to show that $y$ has a neighbor in exactly one $A_{i+2}$-clique. Suppose not. Then there are nonadjacent vertices $a$ and $b$ in $A_{i+2}$ such that $y$ is adjacent to both $a$ and $b$. Then pick a neighbor of $y$ in each  $A_{i-2}$ and $A_{i+1}$, say $p$ and $q$ respectively; but then $\{p,a,q,b,y\}$ induces a $4$-wheel which is a contradiction.  This proves \cref{Ai+2-Ai-2-yi}. $\Diamond$

\smallskip
\no{$(iii)$}:~Let $y\in Y_i$.  We may assume, up to symmetry, that $y$ is not complete to $A_{i-1}$, and let $p$ be a nonneighbor of $y$ in $A_{i-1}$. So by \ref{SA-Yi-Ai-1-Ai+1}, $A_{i-1}$ is a clique.
 Suppose to the contrary that $y$ has a nonneighbor in $A_{i-2}\cup A_{i+2}$, say $q$. If $q\in A_{i-2}$, then for any neighbor of $y$ in $A_{i+1}$, say $r$, we see that $q$-$p$-$v_i$-$r$-$y$ is a $P_5$, a contradiction; so $q\in A_{i+2}$. Pick a neighbor of $y$ in  each $A_{i-1}$ and $A_{i+1}$, say $a$ and $b$ respectively. Since $A_{i-1}$ is a clique, $pa\in E$. Now we see that $p$-$a$-$y$-$b$-$q$ is a $P_5$, a contradiction. This proves the first assertion of \cref{YiYi+2Yi-2comYi+1ancom}, and the second assertion    follows from \cref{Ai+2-Ai-2-yi}. $\Diamond$

\smallskip
\no{$(iv)$}:~Let $y, y' \in Y_i$, and suppose $y,y'$ are nonadjacent. By \ref{yi-Ai-1-Ai+1}, we may assume that $y$ is complete to $A_{i-1}$. Then by the definition of $Y_i$, clearly $y$ and $y'$ have a common neighbor in $A_{i-1}$, say $p$.  So by the definition of $Y_i$ and by \ref{p5}, $y$ and $y'$  have a common neighbor in  $A_{i+2}$, say $q$. By the same argument, if $y$ and $y'$  have a common neighbor in  $A_{i+1}$,
  then they have a common neighbor in $A_{i-2}$.  If $y$ and $y'$ do not share a common neighbor in $A_{i+1}$, then by \cref{YiYi+2Yi-2comYi+1ancom}, $A_{i-2}$ is a clique, and so by \cref{Ai+2-Ai-2-K-yi}, $y$ and $y'$ have a common neighbor in $A_{i-2}$.  In either case, $y$ and $y'$ have a common neighbor in $A_{i-2}$, say $r$. Then $\{p,y,q,y',r\}$ induces a $4$-wheel, a contradiction.  This proves \cref{Yi-clq}. $\Diamond$

\smallskip
\no{$(v)$}:~Let $y\in Y_i$ and $y'\in Y_{i+1}$, and  suppose $y$ and $y'$ are nonadjacent. Let $p$ be a neighbor of $y$ in $A_{i-2}$.  If $py'\notin E$, then for any neighbor of $y'$ in  $A_i$, say $a$,  and for any neighbor of $y$ in $A_{i+1}$, say $b$, $p$-$y$-$b$-$a$-$y'$ is a $P_5$, a contradiction; so we may assume that $py'\in E$.
Also it follows from the definition of $Y_{i+1}$, and by \ref{SA-Yi-Ai-1-Ai+1} and \cref{Ai+2-Ai-2-K-yi}, that $y$ and $y'$ have a common neighbor in $A_{i-1}$, say $q$, and by the same argument, $y$ and $y'$ have a common neighbor in $A_{i+2}$, say $r$. But then $\{y',q,y,r,p\}$ induces a $4$-wheel, a contradiction. So $Y_i$ is complete to $Y_{i+1}$. Likewise, $Y_i$ is complete to $Y_{i-1}$. This proves \cref{YiYi+1-clq}. $\Diamond$

\smallskip
\no{$(vi)$}:~Suppose there is a vertex $t\in T$ which has no neighbor in $X$. Let $Q$ be the vertex-set of the component of $G[T]$ containing $t$. Then by \ref{xi-compT}, $Q$ is anticomplete to $X$.  Then since $G$ is connected,   $N(Q)\cap Y\neq \es$.  Since $N(Q)\cap Y$ is not a clique cutset between $A$ and $Q$, there are nonadjacent vertices, say $y,y' \in N(Q)\cap Y$. Then by \cref{Yi-clq} and \cref{YiYi+1-clq}, we may assume that $y\in Y_2$ and $y'\in Y_5$. Now pick a neighbor of $y$ in $A_5$, say $a$, and a neighbor of $y'$ in $A_2$, say $a'$. But then $a$-$y$-$t$-$y'$-$a'$ is a $P_5$, a contradiction. This proves \cref{F1-TXY2}.
  \hfill{$\Box$}

 \begin{lemma}\label{F1-YiYi+1Aj-2}
For each $i\in [5]$, $Y_i\cup Y_{i+1}$ is complete to exactly one $A_{i-2}$-clique.
\end{lemma}
\no{\it Proof of \cref{F1-YiYi+1Aj-2}}. First we show that for each $i$, $Y_i$ is complete to exactly one $A_{i-2}$-clique. Suppose not.  Then by \cref{Ai+2-Ai-2-K-yi}, \cref{Ai+2-Ai-2-yi} and \cref{Yi-clq}, there are adjacent vertices $y,y'$ in $Y_i$, and nonadjacent vertices $a,b$ in $A_{i-2}$ such that $ya, y'b\in E$ and $yb,y'a\notin E$. Then by \cref{YiYi+2Yi-2comYi+1ancom}, $\{y,y'\}$ is complete to $A_{i+1}$ and $A_{i-1}$. Now if $y$ and $y'$ have a common neighbor in $A_{i+2}$, say $p$, then $\{p,y,v_{i-1},b,y'\}$ induces a $4$-wheel, a contradiction; so we may assume that there is a vertex $q\in A_{i+2}$ such that $yq\in E$ and $y'q\notin E$. But then $\{v_{i+1},q,a,v_{i-1},y',y\}$ induces a $5$-wheel, a contradiction. So for each $i$, $Y_i$ is complete to exactly one $A_{i-2}$-clique.

Now suppose that the lemma is not true. Then  by our preceding argument, there are  $A_{i-2}$-cliques, say $B$ and $D$, such that $B\cap D=\es$, $Y_i$ is complete to $B$, and anticomplete  to $A_{i-2}\sm B$, and $Y_{i+1}$ is complete to  $D$, and anticomplete to $A_{i-2}\sm D$. Then clearly $A_{i-2}$ is not a clique, and so by \cref{YiYi+2Yi-2comYi+1ancom}, $Y_i$ is complete to $A_{i-1}$, and $Y_{i+1}$ is complete to $A_{i+2}$. Now pick a vertex  $y\in Y_{i}$, and a neighbor of $y$ in $A_{i+2}$, say $a$. Also, pick a vertex $y'\in Y_{i+1}$, and neighbor of $y'$ in $A_{i-1}$, say $a'$. But now for any $b\in B$, by \cref{YiYi+1-clq}, $\{y',a,b,a',y\}$ induces a $4$-wheel, a contradiction.  This proves
\cref{F1-YiYi+1Aj-2}. \hfill{$\Box$}

\smallskip
If $Y\neq \es$, by \cref{F1-YiYi+1Aj-2}, for $i\in [5]$,   let $B_{i-2}$ be the $A_{i-2}$-clique such that $Y_i\cup Y_{i+1}$ is complete to $B_{i-2}$, and anticomplete to  $A_{i-2}\sm B_{i-2}$, and   let $B_{i+2}$ be the $A_{i+2}$-clique such that $Y_i\cup Y_{i-1}$ is complete to $B_{i+2}$, and anticomplete to  $A_{i+2}\sm B_{i+2}$.

\begin{lemma} \label{prop-yi}
For $i\in [5]$, the following properties hold:
\begin{lemmalist}
\item\label{YiAi+1Ai-1graded} For $j\in \{i-1,i+1\}$, each $A_j$-clique has a vertex which is complete to $Y_i$.
\item \label{F1-Yi+1-XiXi+2}
$Y_{i+1}$ is anticomplete to $X_{i}\cup X_{i+2}$.
\item \label{Xi-Yi+2Yi-2}
At least one of  $X_i$, $Y_{i+2}\cup Y_{i-2}$ is empty.
\item \label{F1-Yi+1-XiAj}
Each $y\in Y_{i+1}$ and $x\in X_{i}$ have a common neighbor in each $A_i$, $A_{i+2}$ and $A_{i-2}$, and  each $y\in Y_{i+1}$ and
$x\in X_{i+2}$ have a common neighbor in each $A_i$, $A_{i+2}$ and $A_{i-1}$.
\item\label{F1-XiYi+1-Ai}
 If $X_i\neq \es$, then $Y_{i+1}\cup Y_{i-1}$ is complete to $A_i$.
 \item \label{Yi-Yi-2Yi+2-ac}
If $X\neq \es$, then $Y_i$ is anticomplete to $Y_{i+2}\cup Y_{i-2}$.
\item\label{F1-TY2Y5}
If $X\neq \es$, then no vertex in $T$ has neighbors in both $Y_{i-1}$ and $Y_{i+1}$.
\end{lemmalist}
\end{lemma}
\no{\it Proof of \cref{prop-yi}.}~$(i)$:~We prove the statement for $j=i+1$. If $A_{i+1}$ is not a clique, then by \ref{SA-Yi-Ai-1-Ai+1}, $Y_i$ is complete to $A_{i+1}$, and   \cref{YiAi+1Ai-1graded} holds; so assume that  $A_{i+1}$ is a clique. Now if $G[Y_i\cup A_{i+1}]$ contains an induced  $C_4$, say with vertex-set $\{p,q,r,s\}$, then for any  $a\in B_{i+2}$, $\{p,q,r,s,a\}$ induces a $4$-wheel, a contradiction; so $G[Y_i\cup A_{i+1}]$ is $C_4$-free.  Since $Y_i$ is a clique (by \cref{Yi-clq}) and since each vertex in $Y_i$ has a neighbor in $A_{i+1}$ (which is a clique), by Lemma \ref{lem:squarefree}, $A_{i+1}$ has a vertex which is complete to $Y_i$. This proves \cref{YiAi+1Ai-1graded}. $\Diamond$

 \smallskip
\no{$(ii)$}:~Suppose, up to symmetry, there are adjacent vertices, say $y\in Y_{i+1}$ and $x\in X_i$. Pick a neighbor of $y$ in each $A_{i-1}$ and $A_i$, say $p$ and $q$ respectively. If $x$ and $y$ have a  common neighbor in $A_{i-2}$, say $r$,  then, by \ref{XiAi}, $\{q,x,r,p,y\}$ induces a $4$-wheel, a contradiction; so there is a vertex $w\in A_{i-2}$ such that $yw\in E$ and  $xw\notin E$. Then by \ref{xAi+2Ai-2}, $x$ is complete to $A_{i+2}$. Now pick any  neighbor of $y$ in $A_{i+2}$, say $s$. Then, by \ref{XiAi}, $\{p,q,x,s,w,y\}$ induces a $5$-wheel, a contradiction. So $Y_{i+1}$ is anticomplete to $X_i$. Likewise, $Y_{i+1}$ is   anticomplete to $X_{i+2}$. This proves \cref{F1-Yi+1-XiXi+2}. $\Diamond$

 \smallskip
\no{$(iii)$}:~Suppose not. Let $x\in X_i$, and, up to symmetry, let $y\in Y_{i+2}$.   Pick any neighbor of $y$ in $A_{i-1}$, say $p$. It follows from  \ref{Ai+2-Ai-2-K-xi} and \ref{SA-Yi-Ai-1-Ai+1} that $x$ and $y$ have a common neighbor in $A_{i-2}$, say $a$.
  Now if $xy\in E$, then for any neighbor of $y$ in  $A_{i}$, say $a'$, by \ref{XiAi},  $\{p,a,x,a',y\}$ induces a $4$-wheel, a contradiction; so we may assume that $xy\notin E$. Then pick a neighbor of $y$ in $A_{i+1}$, say $b$, and a neighbor of $x$ in  $A_{i+2}$, say $b'$; but then $p$-$y$-$b$-$b'$-$x$ is a $P_5$ which is a contradiction. This proves \cref{Xi-Yi+2Yi-2}.  $\Diamond$

\smallskip
\no{$(iv)$}:~We prove the first assertion, and the proof of the other  is similar. Suppose $y\in Y_{i+1}$ and $x\in X_i$. By \ref{XiAi}, $x$ is complete to $A_i$, and so by the definition of $Y_{i+1}$, $x$ and $y$ have a common neighbor in $A_i$. By \cref{F1-Yi+1-XiXi+2}, we know that $yx\notin E$. Now $x$ and $y$ have a common neighbor in each $A_{i+2}$ and $A_{i-2}$, by \ref{p5}. This proves \cref{F1-Yi+1-XiAj}. $\Diamond$

\smallskip
\no{$(v)$}:~Let $x\in X_i$. Let $y\in Y_{i+1}$ and $a\in A_i$, and suppose $y,a$ are nonadjacent. By \cref{F1-Yi+1-XiXi+2}, $xy\notin E$, and by \cref{F1-Yi+1-XiAj}, $x$ and $y$ have a common neighbor in $A_{i-2}$, say $a'$. Then by \ref{XiAi}, $y$-$a'$-$x$-$a$-$v_{i+1}$ is a $P_5$, a contradiction. So $Y_{i+1}$ is complete to $A_i$. Likewise, $Y_{i-1}$ is complete to $A_i$. This proves   \cref{F1-XiYi+1-Ai}. $\Diamond$

\smallskip
\no{$(vi)$}:~Suppose not. We may assume that there are adjacent vertices, say $y\in Y_i$ and $y'\in Y_{i+2}$. Since $Y_i,Y_{i+2}\neq \es$, by \cref{Xi-Yi+2Yi-2}, $X_{j}=\es$, for $j\neq i+1$. Now we claim that $X_{i+1}=\es$.  Suppose not. Let $x\in X_{i+1}$. Then by \cref{F1-Yi+1-XiXi+2}, $\{y,y'\}$ is anticomplete to $x$, and by \cref{F1-XiYi+1-Ai}, $\{y,y'\}$ is complete to $v_{i+1}$. If $y$ and $y'$ have a common neighbor in $A_{i-2}$, say $a$, then for any neighbor of $y$ in $A_{i+2}$, say $a'$,  $\{a,a',v_{i+1},y',y\}$ induces a $4$-wheel, a contradiction. So we may assume that $y$ and $y'$ do not share a common neighbor in $A_{i-2}$. Now by \cref{F1-Yi+1-XiAj}, $x$ and $y$ have a common neighbor in $A_{i-2}$, say $p$. But then for any neighbor of $y'$ in $A_{i}$, say $q$, we see that $x$-$p$-$y$-$y'$-$q$ is a $P_5$, a contradiction; so $X_{i+1}=\es$. Thus we conclude that $X=\es$, a contradiction to our assumption that $X\neq \es$.  This proves \cref{Yi-Yi-2Yi+2-ac}. $\Diamond$

\smallskip
\no{$(vii)$}:~We prove the assertion for $i=1$.  If some vertex in $T$, say $t$, has neighbors in  both $Y_2$ and $Y_5$, say $y$ and $y'$, respectively. Then by \cref{Yi-Yi-2Yi+2-ac}, $yy'\notin E$. Now pick a neighbor of $y$ in $A_5$, say $a$, and a neighbor of $y'$ in $A_2$, say $a'$, and then $a$-$y$-$t$-$y'$-$a'$ is a $P_5$, a contradiction. This proves \cref{F1-TY2Y5}.  \hfill{$\Box$}

\smallskip
For $i\in [5]$, if $Y_i\cup Y_{i+2}\neq \es$ and if there is a vertex in each $A_{i+1}$-clique which is complete to $Y_i\cup Y_{i+2}$, then we pick one such vertex, and let ${\mathbb A}_{i+1}$ be the union of those vertices; otherwise, we let ${\mathbb A}_{i+1}:=R_{A_{i+1}}$. (In any case, ${\mathbb A}_{i+1}$ is a maximum independent set of $A_{i+1}$.)

\begin{lemma}\label{clq-ABCD}
The set ${\mathbb A}_{i-1}\cup  {\mathbb A}_{i-2}$ meets each maximal clique of $G$ in $G[A_{i-1}\cup A_{i-2}\cup Y_i\cup Y_{i+1}]$ twice. Likewise, ${\mathbb A}_{i+1}\cup  {\mathbb A}_{i+2}$ meets each maximal clique of $G$ in $G[A_{i+1}\cup A_{i+2}\cup Y_i\cup Y_{i-1}]$ twice.
\end{lemma}
\no{\it Proof of \cref{clq-ABCD}}.~By \cref{Yi-clq} and \cref{YiYi+1-clq}, $Y_i\cup Y_{i+1}$ is a clique. Also, we know that $Y_i\cup Y_{i+1}$ is complete to $B_{i-2}$, and anticomplete to  $A_{i-2}\sm B_{i-2}$. Also, $Y_{i+1}$ is complete to $B_{i-1}$, and anticomplete to  $A_{i-1}\sm B_{i-1}$. Let $M$ be a maximal clique in $G[A_{i-2}\cup A_{i-1}\cup Y_i\cup Y_{i+1}]$. If $M$ has no vertex from $Y_i$, clearly the assertion holds. So $M\cap Y_i\neq \es$. If $M$ has no vertex from $Y_{i+1}$, then  $M$ is of the form $Y_i\cup B_{i-2}\cup D_{i-1}$, where $D_{i-1}$ is a subset of some $A_{i-1}$-clique $A^*$, and is the set of vertices in $A^*$ which are complete to $Y_i$ (by \cref{YiAi+1Ai-1graded}).  Since ${\mathbb A}_{i-2}\cup  {\mathbb A}_{i-1}$ contains vertices from both $B_{i-2}$ and $D_{i-1}$, the claim holds.  Finally, if $M\cap Y_{i+1}\neq \es$, then  by \cref{YiYi+1-clq}, $M$ is of the form $Y_i\cup Y_{i+1}\cup B_{i-2}\cup D_{i-1}$, where $D_{i-1}$ is a subset of  $B_{i-1}$, and is the set of vertices in $B_{i-1}$ which are complete to $Y_i$ (by \cref{YiAi+1Ai-1graded}). So again, as earlier, ${\mathbb A}_{i-2}\cup  {\mathbb A}_{i-1}$ meets $M$ twice. This proves \cref{clq-ABCD}. \hfill{$\Box$}

\begin{lemma}\label{lem:C5-YneXe}
If $Y$ is nonempty, and $X$ is empty, then $G$ is  either a nice graph or a quasi-line graph.
\end{lemma}
\no{\it Proof of \cref{lem:C5-YneXe}}.~
Since $X$ is empty,  by \cref{F1-TXY2}, $T=\es$. Now:

\begin{claim}\label{YiYi+2-Ai+1-clq}
If there is an $i\in [5]$ such that $Y_i$ and $Y_{i+2}$ are not complete to $A_{i+1}$, then  $A_i$ is a clique, for all $i\in [5]$.
 \end{claim}
\no{\it Proof of \cref{YiYi+2-Ai+1-clq}}.~Since  $Y_i$ is not complete to $A_{i+1}$, by \cref{SA-Yi-Ai-1-Ai+1}, $A_{i+1}$ is  a clique, and by \cref{YiYi+2Yi-2comYi+1ancom}, $A_{i+2}\cup A_{i-2}$ is a clique. Likewise, since $Y_{i+2}$ is not complete to $A_{i+1}$, by \cref{YiYi+2Yi-2comYi+1ancom}, $A_{i-1}\cup A_{i}$ is a clique. Thus we conclude that $A_i$ is a clique, for all $i$. This proves \cref{YiYi+2-Ai+1-clq}. $\Diamond$

 \begin{claim}\label{AiClq-3K1-free}
 If $A_i$ is a clique, for all $i\in [5]$, then $G$ is $3K_1$-free.
 \end{claim}
 \no{\it Proof of \cref{AiClq-3K1-free}}.~Suppose  that $G$ contains a triad, say $\{u,v,w\}$. Since $G[A]$ is $3K_1$-free, we may assume that $u\in Y_j$, for some $j$. Then by \ref{yi-Ai-1-Ai+1} and \cref{Ai+2-Ai-2-K-yi},  $u$ is complete to either $A_{j+1}\cup A_{j+2}\cup A_{j-2}$ or $A_{j+2}\cup A_{j-1}\cup A_{j-2}$; we may assume, without loss of generality,  that $u$ is complete to $A_{j+1}\cup A_{j+2}\cup A_{j-2}$. Since  $A_{j}\cup A_{j-1}\cup Y_{j+2}$ is a clique (by \cref{Ai+2-Ai-2-K-yi}), and since $Y_j$ is complete to
$Y_{j+1}\cup Y_{j-1}$ (by \cref{YiYi+1-clq}), one of $v,w$ belongs to $Y_{j-2}$; and we may assume that  $v\in Y_{j-2}$. Then by  \cref{Ai+2-Ai-2-K-yi},   $v$ is complete to $A_j\cup A_{j+1}$. So $w\in A_{j-1}$. But then  for any neighbor of $u$ in $A_{j-2}$, say $a$, and for any neighbor of $v$ in $A_j$, say $b$, we see that $u$-$a$-$w$-$b$-$v$ is a $P_5$, a contradiction. This proves \cref{AiClq-3K1-free}. $\Diamond$

\smallskip
First suppose that there is an $i\in [5]$ such that $Y_i$ is not anticomplete to $Y_{i+2}$. Let $y\in Y_i$ and $y'\in Y_{i+2}$ be adjacent.
Suppose $y$ and $y'$ share a common neighbor in $A_{i+1}$, say $a$. We know, by \ref{SA-Yi-Ai-1-Ai+1} and \cref{Ai+2-Ai-2-K-yi}, that $y$ and $y'$ share a common neighbor in $A_{i-1}$, say $a'$. Then for a neighbor of $y'$ in $A_i$, say $a''$,   $\{a,y,a',a'',y'\}$ induces a $4$-wheel, a contradiction; so suppose that $y$ and $y'$ do not share a common neighbor in $A_{i+1}$. So $y$ and $y'$ are not complete to $A_{i+1}$, hence $Y_i$ and $Y_{i+2}$ are not complete to $A_{i+1}$. Then by \cref{YiYi+2-Ai+1-clq}, $A_i$ is a clique for all $i\in [5]$, and then, by \cref{AiClq-3K1-free}, $G$ is $3K_1$-free. So, by  \cref{3K1W4-free-Struc}, $G$ is either a quasi-line graph or  a nice graph, and we are done.

Next we may assume that for each $i\in [5]$, $Y_i$ is anticomplete to $Y_{i+2}\cup Y_{i-2}$. By \cref{YiYi+1-clq}, $Y_i$ is complete to $Y_{i+1}\cup Y_{i-1}$. See Figure~\ref{fig:xory}:$(c)$ for a sketch of $G$. Also,  we may assume that if $Y_i$ and $Y_{i+2}$ are nonempty, then at least one of $Y_i$, $Y_{i+2}$ is complete to $A_{i+1}$  (for, otherwise, by \cref{YiYi+2-Ai+1-clq} and \cref{AiClq-3K1-free}, $G$ is $3K_1$-free, and we conclude using \cref{3K1W4-free-Struc}).
Now we define three sets $S_1:={\mathbb A}_1\cup {\mathbb A}_3$, $S_2:={\mathbb A}_2\cup {\mathbb A}_4$ and $S_3:={\mathbb A}_5$. Then $S_1,S_2$ and $S_3$ are stable sets. Clearly, for $i\in [5]$,   by \cref{Yi-clq} and \cref{F1-YiYi+1Aj-2}, ${\mathbb A}_{i+2}\cup {\mathbb A}_{i-2}$ meets each maximal clique of $G$ in $G[Y_i\cup A_{i+2}\cup A_{i-2}]$   twice, and  by \cref{clq-ABCD}, ${\mathbb A}_{i+1}\cup  {\mathbb A}_{i+2}$ meets each maximal clique of $G$ in $G[Y_i\cup Y_{i-1}\cup A_{i+1}\cup A_{i+2}]$  twice.  So we conclude that $S_1\cup S_2\cup S_3$ meets each maximum clique of $G$ at least twice, and other maximal cliques at least once, and that $G$ is nice.
This completes the proof of \cref{lem:C5-YneXe}.  \hfill{$\Box$}

\smallskip
Now we  prove our main result of this subsection, and is given below.

\begin{thm}\label{thm:P5wheel-free-C5}
If a connected ($P_5$,\,wheel)-free atom  $G$  contains an induced $C_5$, then  $G$ is either a nice graph or   a quasi-line graph.
\end{thm}
\no{\it Proof of Theorem~\ref{thm:P5wheel-free-C5}}.~Let $G$ be a connected ($P_5$,\,wheel)-free atom that   contains an induced $C_5$,
 say  $v_1$-$v_2$-$v_3$-$v_4$-$v_5$-$v_1$.     Then we define the sets $A$, $X$, $Y$, $Z$ and $T$ as in Section~\ref{sec-C5} with $v_i\in A_i$, for each $i$, and we use the facts \ref{part}--\ref{TXi-nonadj-last} shown in Section~\ref{sec-C5}, and properties in  \cref{prop-xi,lem:XiAi+2maxclqyemp,prop-yi1,F1-YiYi+1Aj-2,prop-yi,clq-ABCD}.
Since $Z=\es$,   if $T\neq \es$, then   $G[T]$ is $P_3$-free (by \ref{5WF-T-P3free}).
 Let $\cal M$ be the set of   maximum cliques of $G$.   Let $L$ consist of one vertex from each $T$-clique, otherwise let $L:=\es$, and let $L'$ consist of one vertex (which is not in $L$)  from each nontrivial $T$-clique, otherwise let $L':=\es$. By \ref{G7-clm1}, for $j\in [5]$, if $X_j, X_{j+1}\neq \es$,  $R_{X_j}\cup R_{X_{j+1}}\cup R_{A_{j-2}}$ meets each maximal clique of $G$ in $G[X_j\cup X_{j+1}\cup A_{j-2}]$ at least twice.
   Now if $X \cup Y=\es$, then since $G$ is connected, by \cref{F1-TXY2},  $T=\es$, and then the sets $S_1:= R_{A_1}\cup R_{A_3}$, $S_2:=R_{A_2}\cup R_{A_4}$ and $S_3:=R_{A_5}$  are the desired stable sets, and we are done.
If one of $X$, $Y$ is empty, then the theorem follows from \cref{lem:C5-YneXe,lem:C5-Xne-Ye}. So we may assume that both $X$ and $Y$ are  nonempty.
Now if $Y_{i+1}\neq \es$, then by \cref{F1-Yi+1-XiAj} and \ref{Ai+2-Ai-2-K-xi},  $X_{i}$ is complete to $B_{i-2}$, and $X_{i+2}$ is complete to $B_{i-1}$. Recall that since $X\neq\es$, by \cref{Yi-Yi-2Yi+2-ac}, $Y_i$ is anticomplete to $Y_{i+2}\cup Y_{i-2}$. Now we split the proof   into two cases.

\medskip
\no{\bf Case~1}~\label{oneXYe} For each $i\in [5]$, one of $X_i$, $Y_i$ is empty.

\smallskip
Since $Y\neq \es$, let  $Y_2\neq \es$; so  $X_2=\es$. By \cref{Xi-Yi+2Yi-2}, $X_4\cup X_5=\es$. Since $X\neq \es$, $X_1\cup X_3\neq \es$; we may assume that $X_1\neq\es$; so $Y_1=\es$. Again by \cref{Xi-Yi+2Yi-2}, $Y_3\cup Y_4=\es$.
By \cref{F1-XiYi+1-Ai}, $Y_2\cup Y_5$ is complete to $A_1$.
By \cref{F1-Yi+1-XiXi+2}, $Y_{2}\cup Y_5$ is anticomplete to $X_1\cup X_3$.  Recall that  $Y_2$ is complete to $B_4\cup B_5$, and anticomplete to $(A_4\sm B_4)\cup (A_5\sm B_5)$, and since $Y_2\neq\es$, $X_1$ is complete to $B_4$.    Moreover, we have the following:

\begin{cla}\label{2.1:T-neq-emp}
If $T\neq \es$, then the following hold: (a)  $T$ is complete to $Y_2$.  (b) For $j\in \{1,3\}$, given an $X_j$-clique, $X_j^*$, each $T$-clique  is either complete or anticomplete to $X_j^*$. (c) $Y_5=\es$.
 \end{cla}
\no{\it Proof of \cref{2.1:T-neq-emp}}.~$(a)$:~Let $T'$ be a $T$-clique in $G$. Then by \cref{F1-TXY2}, $N(T')\cap (X_1\cup X_3)\neq \es$. Since $Y_2$ is anticomplete to $X_1\cup X_3$ (by \cref{F1-Yi+1-XiXi+2}), it follows from \ref{TXi-nonadj-last}:(a) that $T'$ is complete to $Y_2$. This proves $(a)$, since $T'$ is arbitrary.

\no $(b)$:~Since $Y_2$ is anticomplete to $X_1\cup X_3$, $(b)$ follows from $(a)$,  \cref{F1-TXY2}, and \ref{TXi-nonadj-last}:(b).

\no $(c)$:~Suppose that $Y_5\neq \es$. Then, by \cref{Xi-Yi+2Yi-2}, $X_3=\es$. Let $t\in T$. Then by \cref{F1-TXY2}, $t$ has a neighbor in $X_1$. Since $Y_2\cup Y_5$ is anticomplete to $X_1$ (by \cref{F1-Yi+1-XiXi+2}), it follows from  \ref{TXi-nonadj-last}  that $t$ has neighbors in both $Y_2$ and $Y_5$ which is a contradiction to \cref{F1-TY2Y5}. This proves $(c)$.  $\Diamond$

\begin{cla}\label{2.1:Y2T-A4A5}
 If $T\neq \es$, then $Y_2$ is complete to either $A_4$ or $A_5$. So, if $T\neq \es$, either $A_4$ or $A_5$ is a clique.
\end{cla}
\no{\it Proof of \cref{2.1:Y2T-A4A5}}.~Suppose not. Then  there are vertices $y\in Y_2$, $p\in A_4$, and $q\in A_5$ such that $yp,yq\notin E$. Let $t\in T$. Then by \cref{2.1:T-neq-emp}:(a),   $yt\in E$. But then since $Y_2$ is complete to $A_1$, for any neighbor of $y$ in $A_1$, say $r$, we see that $p$-$q$-$r$-$y$-$t$ is a $P_5$, a contradiction. This proves \cref{2.1:Y2T-A4A5}.   $\Diamond$

\begin{cla}\label{X1A3-X3A1-11clq}
If $K$ is an $X_1$-clique, and  $D$ is an $A_3$-clique, then either $K$ is complete to $D$ or $K$ is anticomplete to $D$.
Likewise, if $K'$ is an $X_3$-clique, and  $D'$ is an $A_1$-clique, then either $K'$ is complete to $D'$ or $K'$ is anticomplete to $D'$.
  \end{cla}
 \no{\it Proof of \cref{X1A3-X3A1-11clq}}.~Suppose that  $K$ is not anticomplete to $D$. Then, there is an  $x\in K$ which has a neighbor in $D$. Let $a\in D$ be such that $a\in {\mathbb A}_3$ (such a vertex exists, by \cref{YiAi+1Ai-1graded}). Then by \ref{Ai+2-Ai-2-K-xi}, $x$ is complete to $D$; so $xa\in E$.  Let  $x'(\neq x)\in K$ be arbitrary. We claim that $x'$ is complete to $D$. Suppose not.  Then again by \ref{Ai+2-Ai-2-K-xi}, $x'$ is anticomplete to $D$; so $x'a\notin E$. But then, for any $y\in Y_2$ and $b\in B_5$, $x'$-$x$-$a$-$y$-$b$ is a $P_5$, a contradiction. So $x'$ is complete to $D$. Since $x'$ is arbitrary, $K$ is complete to $D$. This proves  \cref{X1A3-X3A1-11clq}.  $\Diamond$

 \begin{cla}\label{X1A4-X3A5-11clq}
	If $K$ is an $X_1$-clique, and  $D$ is an $A_4$-clique, then either $K$ is complete to $D$ or $K$ is anticomplete to $D$.
	Likewise, if $K'$ is an $X_3$-clique, and  $D'$ is an $A_5$-clique, then either $K'$ is complete to $D'$ or $K'$ is anticomplete to $D'$.
\end{cla}
\no{\it Proof of \cref{X1A4-X3A5-11clq}}.~Suppose that  $K$ is not anticomplete to $D$. We may assume that $D\neq B_4$. Then, there is an  $x\in K$ which has a neighbor in $D$, say $a$. Let  $x'(\neq x)\in K$ be arbitrary. We claim that $x'$ is complete to $D$. Suppose not.  Then, by \ref{Ai+2-Ai-2-K-xi}, $x'$ is anticomplete to $D$; so $x'a\notin E$. But then, for any $y\in Y_2$ and $b\in B_5$, $x'$-$x$-$a$-$b$-$y$ is a $P_5$, a contradiction. So $x'$ is complete to $D$. Since $x'$ is arbitrary, $K$ is complete to $D$. This proves  \cref{X1A4-X3A5-11clq}. $\Diamond$

\smallskip
Now consider any maximum clique of $G$ in $G[X_i\cup A_{i+2}\cup A_{i-2}]$, say $M$. Then by \cref{X1A3-X3A1-11clq} and  \cref{X1A4-X3A5-11clq}, $M\cap A_{i+2}, M\cap A_{i-2}\neq \es$, $M\cap A_{i+2}$ is an $A_{i+2}$-clique and $M\cap A_{i-2}$ is an $A_{i-2}$-clique.

\begin{figure}[h]
\centering
 \includegraphics[height=3.5cm, width=15cm]{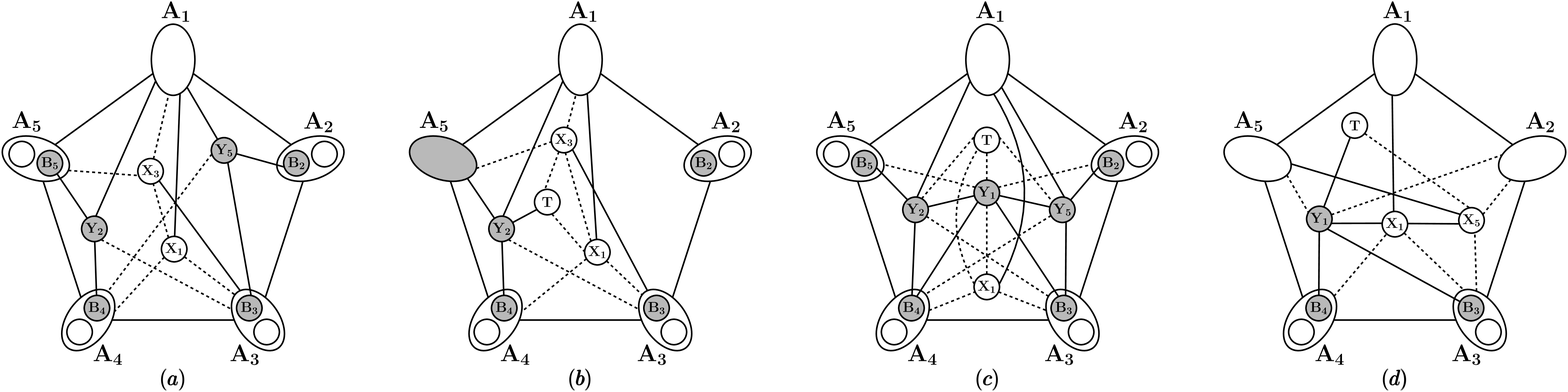}
\caption{Sketch of the graph $G$ in Theorem~\ref{thm:P5wheel-free-C5}: $(a)$ Case~1 when $T=\es$. $(b)$  Case~1 when $T\neq \es$. $(c)$~Case~2 when $X_2\cup X_5=\es$. $(d)$   Case~2 when $X_5\neq \es$.} \label{fig:yx}
\end{figure}

\smallskip

See Figure~\ref{fig:yx}:$(a)$ and Figure~\ref{fig:yx}:$(b)$. By \cref{XiXi+2-11clq} and \cref{X1A3-X3A1-11clq}, we conclude that each $X_1$-clique is either complete or anticomplete to $K$, where is $K$ is an $A_3$-clique or an $X_3$-clique. Likewise, each $X_3$-clique is either complete or anticomplete to $K'$, where  $K'$ is an $A_1$-clique or an $X_1$-clique.  Thus, by \ref{XiAi} and \ref{new-obs}, for $j\in \{1,3\}$, $R_{X_1}\cup R_{X_3}\cup R_{A_j}$ meets each maximal clique of $G$ in $G[X_1\cup X_3\cup A_j]$ at least twice.  Also, by \cref{XiXi+2-11clq}, \cref{2.1:T-neq-emp}:(b), and by  \ref{new-obs}, $R_{X_1}\cup R_{X_3}\cup R_T$ meets each maximal clique of $G$ in $G[X_1\cup X_3\cup  T]$ at least twice. Clearly   $R_{A_1}\cup R_{A_5}$ meets each maximal clique of $G$ in $G[Y_2\cup A_1\cup A_5]$  twice. Likewise,  $R_{A_1}\cup R_{A_2}$ meets each maximal clique of $G$ in $G[Y_5\cup A_1\cup A_2]$    twice.

 Now if $T=\es$, then  using \cref{Xi-O3} and by \cref{clq-ABCD}, we see that the sets $S_1:=R_{A_2}\cup R_{A_5} \cup R_{X_1}$, $S_2:=R_{A_1}\cup {\mathbb A}_3$ and $S_3:={\mathbb A}_4\cup R_{X_3}$ are the desired stable sets. So we may assume that $T\neq \es$. By \cref{2.1:T-neq-emp}, $Y_5=\es$.
By \cref{2.1:Y2T-A4A5}, up to relabelling, we may assume that $A_5$ is a clique. Now we let $S_1:=R_{A_2}\cup R_{A_5} \cup R_{X_1}$, $S_2:= R_{A_1}\cup {\mathbb A}_3\cup R_{T}$ and $S_3:=R_{A_4\sm B_4}\cup R_{X_3}\cup R_{Y_2}$. Clearly, $S_1,S_2$ and $S_3$ are stable sets. Let $S:=S_1\cup S_2\cup S_3$. To justify that $S$ meets each maximum clique of $G$ that has vertices from $A_3\cup B_4$, at least twice, we need the following.

\begin{cla}\label{B4A3-Y2X1}
 If $M \in \cal M$ has vertices from both $B_4$ and  an $A_3$-clique $A_3^*$, and no vertex from $Y_2$,  then $M$ is of the form $X_1^*\cup B_4\cup A_3^*$, where $X_1^*$ is an $X_1$-clique, and $B_4\cup A_3^*\notin {\cal M}$.
\end{cla}
\no{\it Proof of \cref{B4A3-Y2X1}}.~If $X_1^*$ is complete to $A_3^*$, then $M$ is of the form $X_1^*\cup B_4\cup A_3^*$, and we are done. So, we may assume, by \cref{X1A3-X3A1-11clq}, that $X_1^*$ is anticomplete to $A_3^*$. Then, by \ref{Ai+2-Ai-2-K-xi}, $A_3$ is not a clique, and then by \ref{SA-Yi-Ai-1-Ai+1},   $A_3$ is complete to $Y_2$. But now   since $Y_2\cup B_4\cup A_3^*$ is a larger clique than $B_4\cup A_3^*$, we have $B_4\cup A_3^*\notin {\cal M}$.   This proves \cref{B4A3-Y2X1}.  $\Diamond$

\smallskip
 Now by \cref{Xi-O3} and \cref{B4A3-Y2X1}, $R_{X_1}\cup R_{A_4\sm B_4}\cup R_{A_3}$ meets each maximum clique of $G$ in $G[X_1\cup A_3\cup A_4]$ twice. Clearly, $R_{A_5}\cup R_{Y_2}$ meets each maximum  clique of $G$ in $G[Y_2\cup B_4\cup A_5]$ twice, ${\mathbb A}_3\cup R_{Y_2}$ meets each maximum clique of $G$ in $G[Y_2\cup A_3\cup A_4]$ twice (by  \cref{clq-ABCD}), $R_{A_5}\cup R_{A_4\sm B_4}$ meets each maximum clique of $G$ in $G[(A_4\sm B_4)\cup A_5]$, and $R_{Y_2}\cup R_{T}$ meets each maximum  clique of $G$ in $G[Y_2\cup T]$ twice.
Thus, we conclude  that $S$  meets each maximum clique of $G$ at least twice, and other maximal cliques at least once, and that $G$ is nice.

\medskip

\no{\bf Case~2}~\label{lastcase}
There is an index $i\in [5]$ such that $X_i$ and $Y_i$ are nonempty.

\smallskip
 Let $i=1$.  Then by \cref{Xi-Yi+2Yi-2}, $X_{3}\cup X_{4}\cup Y_{3}\cup Y_{4} =\es$.   Recall that $Y_1$ is anticomplete to $X_2\cup X_5$ (by \cref{F1-Yi+1-XiXi+2}), and complete to $Y_2\cup Y_5$ (by \cref{YiYi+1-clq}). Also $X_1$ is complete to $X_2\cup X_5$ (by \ref{XiXi+1}). By \cref{F1-YiYi+1Aj-2}, $Y_1\cup Y_2$ is complete to $B_4$, and anticomplete to $A_4\sm B_4$; $Y_1\cup Y_5$ is complete to $B_3$, and anticomplete to $A_3\sm B_3$.  Also, by \cref{F1-Yi+1-XiAj} and \ref{Ai+2-Ai-2-K-xi}, $X_2$ is complete to $B_4$,  and $X_5$ is complete to $B_3$. Note that since $Y_1\neq \es$, $B_3, B_4\neq \es$. Since $X_1\neq \es$, by \cref{F1-XiYi+1-Ai}, $Y_2\cup Y_5$ is complete to $A_1$. Recall that since $Z=\es$, each vertex in $X\cup Y$
is either anticomplete or good with respect to $T$. By \cref{lemmaXiAi+2meettwice}, $R_{X_i}\cup R_{A_{i+2}}\cup R_{A_{i-2}}$ meet each maximum clique of $G$ in $G[X_i\cup A_{i+2}\cup A_{i-2}]$  twice, and other maximal cliques in $G[X_i\cup A_{i+2}\cup A_{i-2}]$  once.  To proceed further we claim following:
\begin{cla}\label{F1-xiYicom}
	Suppose $x\in X_1$ has a neighbor in  $(A_3\sm B_{3})\cup (A_4\sm B_{4})$. Then $x$ is complete to $Y_1$.
\end{cla}
\no{\it Proof of \cref{F1-xiYicom}}.~We may assume, up to symmetry, that $x$ has a neighbor in $A_3\sm B_3$, say $p$. Let $y\in Y_1$, and suppose $x,y$ are nonadjacent. Now pick a neighbor of $y$ in $A_5$, say $a$.
Then for any $a'\in A_1$, by \ref{XiAi}, $p$-$x$-$a'$-$a$-$y$ is a $P_5$, a contradiction.  This proves \cref{F1-xiYicom}. $\Diamond$

\begin{cla}\label{XYmaxclq}
Let $M$ be a maximal clique of $G$ containing at least one vertex from each of $X_1$ and $Y_1$, and no vertex from $T$. Then $R_{B_3}\cup R_{B_4}\cup R_{X_1}$ meets $M$ at least twice.
\end{cla}
 \no{\it Proof of \cref{XYmaxclq}}.~Let $M\cap X_1=X_1^*$ and let $D$ be the $X_1$-clique such that $X_1^*\subseteq D$. Recall that $Y_1$ is complete to $B_3\cup B_4$. Now we claim that $D$ is complete to  either $B_3$ or  $B_4$. Suppose not.  Then by \ref{x-good} and \ref{xAi+2Ai-2}, there are vertices $x,x'\in D$ such that $x$ is anticomplete to $B_3$, and $x'$ is anticomplete to $B_4$.  Then by the definition of $X_1$, $x$ has a neighbor in $A_3\sm B_3$, and  $x'$ has a neighbor in $A_4\sm B_4$. So, by \cref{F1-xiYicom}, $\{x,x'\}$ is complete to $Y_1$. Then by \cref{F1-xiYicom} and \ref{xAi+2Ai-2}, for any $y\in Y_1$, $b\in B_3$ and $b'\in B_4$ then $\{x,b',b,x',y\}$ induces a $4$-wheel, a contradiction; so $D$ is  complete to either $B_3$ or  $B_4$.  Now since $M$ is a maximal clique and $X_1^*$ is either complete to $B_3$ or $B_4$, we conclude that $M\cap (B_3\cup B_4)\neq \es$. If $X_1^*$ is complete to $B_3\cup B_4$, then clearly the assertion holds. So we assume that $X_1^*$ is not complete to $B_4$, then there is an $x\in X_1^*$  such that  $x$ is anticomplete to $B_4$ (by \ref{Ai+2-Ai-2-K-xi}). So by a previous argument, $D$ is complete to $B_3$. Next we claim that $D$ is complete $Y_1$. Suppose there are vertices $x'\in D$ and $y\in Y_1$ such that $x'y\notin E$. By \cref{F1-xiYicom}, $x\neq x'$, and $x'$ is anticomplete to $A_4\sm B_4$, so by the definition of $X_1$, $x'$ must have a neighbor in $B_4$. Then by \cref{F1-xiYicom} and \ref{x-good}, for any $a\in B_3$ and $a'\in B_4$, $\{a',x',x,y,a\}$ induces a $4$-wheel, a contradiction. So $D$ is complete $Y_1$. Since $B_3\cup D\cup Y_1$ is a clique, and $M$ is a maximal clique, we have $X_1^*=D$,  and hence $M=B_3\cup D\cup Y_1$. Then clearly   $R_{B_3}\cup R_{B_4}\cup R_{X_1}$ meets $M$ at least twice. This proves \cref{XYmaxclq}. $\Diamond$

\begin{cla}\label{Qqeg2-1}
Suppose that $X_2\cup X_5=\es$, and let $Q$ be a $T$-clique.   If there is an $M\in \cal{M}$ such that $Q\subseteq M$, then $|Q|\geq 2$.
\end{cla}
\no{\it Proof of \cref{Qqeg2-1}}.~Recall that each vertex in $X\cup Y$ is either complete or anticomplete to $Q$. For $j\in \{1,2,5\}$, let $Y_j^*:= N(Q)\cap Y_j$. First suppose that $(M\sm Q)\cap X=\es$. If $M\sm Q =Y_1^*\cup Y_2^*$, then $|Y_1^*\cup Y_2^*\cup Q|\geq |B_4\cup Y_1\cup Y_2\cup \{b_3\}|$, where $b_3\in B_3$ is the vertex which is complete to $Y_2$ (by \cref{YiAi+1Ai-1graded}), hence $|Q|\geq 2$. Likewise, if $M\sm Q =Y_1^*\cup Y_5^*$, then $|Q|\geq 2$. So we assume $(M\sm Q)\cap X\neq \es$. Then $M\sm Q=X_1^*\cup Y_1^*$ where $X_1^*$ is a subset of some $X_1$-clique such that $N(Q)\cap X_1^*\neq \es$. Then by \ref{xi-pure}, $X_1^*$ is complete to $A_3\cup A_4$; in particular $X_1^*$ is complete to $B_3\cup B_4$. So  $|X_1^*\cup Y_1^*\cup Q|\geq |B_3\cup B_4\cup X_1^*\cup Y_1^*|$, thus $|Q|\geq 2$.  This proves \cref{Qqeg2-1}.  $\Diamond$

\medskip
First suppose that $X_2\cup X_5=\es$, and we apply \cref{Qqeg2-1}. We refer to Figure~\ref{fig:yx}:$(c)$ for a sketch of the graph $G$.  Then, by \cref{clq-ABCD}, $\cup_{i=1}^5 {\mathbb A}_i$ meets each  maximum clique of $G$ in $G[A\cup Y]$  twice. Since $X_1$ is anticomplete to $Y_2\cup Y_5$, by \cref{XYmaxclq},  and by \cref{clq-ABCD}, ${\mathbb A}_3\cup {\mathbb A}_4\cup R_{X_1}$ meets each maximum clique of $G$ in $G[X\cup Y\cup A_3\cup A_4]$ twice, and clearly, by \ref{XiAi}, ${\mathbb A}_1\cup R_{X_1}$ meets each maximum clique of $G$ in $G[X_1\cup A_1]$ twice. Then since each vertex in $X\cup Y$ is either complete or anticomplete to a $T$-clique, by \cref{Qqeg2-1}, we conclude that $S_1:={\mathbb A}_2\cup {\mathbb A}_5\cup R_{X_1} $, $S_2:={\mathbb A}_1\cup {\mathbb A}_4\cup L$ and $S_3:={\mathbb A}_3\cup L'$ are the required stable sets.  Next suppose that $X_2\cup X_5\neq \es$. We may assume, up to symmetry, that $X_5\neq \es$. Then, by \cref{Xi-Yi+2Yi-2}, $Y_2=\es$.  Next we claim that $X_1$ is complete to $Y_1$. Suppose to the contrary that there are nonadjacent vertices, say $x\in X_1$ and $y\in Y_1$. Then,  by \cref{F1-xiYicom}, $x$ is complete to $B_3\cup B_4$. Now pick any  $x'\in X_5$,  $a\in B_3$, $a'\in B_4$, and pick a common neighbor of $x'$ and $y$ in $A_2$, say $a''$ (by \cref{F1-Yi+1-XiAj}). Then since $X_5$ is complete to $B_3$, we see that $\{y, a', x, x', a'', a\}$ induces a $5$-wheel, a contradiction; so $X_1$ is complete to $Y_1$.
Further, if there are adjacent vertices, say $x\in X_1$ and $b\in B_3$, then for any $x'\in X_5$, $y\in Y_1$, by \cref{F1-Yi+1-XiAj}, $x'$ and $y$ have a common neighbor in $A_2$, say $a$, and then, by \ref{XiXi+1}, $\{x,y,a,x',b\}$ induces a $4$-wheel, a contradiction; so $X_1$ is anticomplete to $B_3$. Likewise, if $X_2\neq \es$, then $X_1$ is anticomplete to $B_4$, a contradiction to \ref{xAi+2Ai-2}; so $X_2=\es$. Since $X_1$ is anticomplete to $B_3$, by \cref{F1-Yi+1-XiAj}, $Y_5=\es$, and  by \ref{xi-pure}, $X_1$ is anticomplete to $T$. By \cref{F1-TXY2}, each vertex in $T$ has a neighbor in $X_5$, and so by \ref{TXi-nonadj-last},  $T$ is complete to $Y_1$. Hence again by \ref{TXi-nonadj-last}, each $T$-clique is either complete or anticomplete to an $X_5$-clique. See Figure~\ref{fig:yx}:$(d)$ for a sketch of the graph $G$.  Moreover, if there is a $M\in \cal M$ which has vertices from a $T$-clique  $T^*$ and from $Y_1$, then  $|Y_1\cup T^*|\geq |Y_1\cup B_3\cup B_4|$, and so $|T^*|\geq 2$.
Now we let $S_1:={\mathbb A}_2\cup {\mathbb A}_5\cup R_{X_1}\cup L $, $S_2:=R_{A_1}\cup R_{A_4}\cup R_{X_5}$ and $S_3:=R_{A_3}\cup L'$.
Then $R_{X_5}\cup L\cup L'$ meets each maximum clique of $G$ in $G[X\cup Y\cup T]$ at least twice, and the other maximal cliques  once.  Also, by \cref{XYmaxclq}, $R_{A_3}\cup R_{A_4}\cup R_{X_1}$ meets each maximum clique of $G$ in $G[X_1\cup Y_1\cup A_3\cup A_4]$ twice, and the other maximal cliques  once. By \ref{G7-clm1}, $R_{X_1}\cup R_{X_5}\cup R_{A_3}$ meets each maximal clique of $G$ in $G[X_1\cup X_5\cup A_3]$ at least twice.
Now by using  \cref{clq-ABCD}, we observe that $S_1\cup S_2\cup S_3$  meets each maximum clique of $G$ at least twice, and other maximal cliques at least once. So  $G$ is nice.   This completes the proof of Theorem~\ref{thm:P5wheel-free-C5}. \hfill{$\Box$}

\subsection{Structure of ($P_5$,\,$C_5$,\,$4$-wheel)-free  graphs that contain an induced $C_7^c$}
Let $C^*$ be the $C_7^c$ with vertices $v_1,v_2,\ldots,v_7$  and edges $v_iv_{i+1}$ and  $v_iv_{i+2}$ for each $i$ modulo $7$. Let $H^*$ be the graph obtained from $C^*$ by adding two vertices $v_8$ and $v_9$ and edges $v_8v_1,v_8v_2,v_8v_5, v_9v_5,v_9v_6$ and $v_9v_2$.

\begin{thm}\label{thm-C7c}
If a connected ($P_5$,\,$C_5$,\,$4$-wheel)-free graph $G$ contains an induced $C_7^c$, then $G$ is a blowup of $H^*$, and hence $G$ is nice.
\end{thm}
\no{\it Proof of Theorem~\ref{thm-C7c}}.~For convenience, we consider the complement graph of $G$, say $H$. So $H$ is a ($P_5^c$,\,$C_5$,\,$2K_2\cup K_1$)-free graph such that $H^c (\cong G)$ is connected, and contains an induced  $C_7$, say $u_1$-$u_2$-$u_3$-$u_4$-$u_5$-$u_6$-$u_7$-$u_1$.  So we may assume that there are seven nonempty and pairwise disjoint sets $A_1, ..., A_7$ such that for each $i$ modulo~$7$ the
set $A_i$ is complete to $A_{i-1}\cup A_{i+1}$, and anticomplete to
$A_{i-2}\cup A_{i-3}\cup A_{i+2}\cup A_{i+3}$.  Let $A:=A_1\cup\cdots\cup A_7$.  We choose
these sets such that $A$ is maximal, and let $u_i\in A_i$.  For each $i\in [7]$, let $B_i$ denote the set $\{x\in V(H)\sm A \mid x$  has a neighbor in each $A_j$, $j\in \{i, i+1,i+2,i+3\}$, and  $x$  is anticomplete to $A_{i-1}\cup A_{i-2}\cup A_{i-3}\}$.
Let $B:=B_1\cup\cdots\cup B_7$. Let $D$ denote the set $\{x\in V(H)\sm A \mid x \mbox{ has a neighbor in } A_i, \mbox{ for each } i\in [7]\}$. Clearly, since the graph $H$ is ($P_5^c$,\,$C_5$)-free, we have the following simple observation:

 \begin{pro} \label{c5-free}
Let $P$ be a $P_4$ in $H$, say $a_1$-$a_2$-$a_3$-$a_4$. Then any vertex in $V(H)\sm V(P)$ which is adjacent to both
$a_1$ and $a_4$, is adjacent to both $a_2$ and $a_3$.
\end{pro}

Moreover, the following hold, for each $i\in [7]$:

 \begin{pro}\label{C7-neiA}
{\it Each vertex in $V(H)\sm A$ has a neighbor in $A$.}
\end{pro}
\no{\it Proof of $\ref{C7-neiA}$}.~If some $x\in V(H)\sm A$ has no neighbor in $A$, then $\{u_1,u_2,u_4,u_5,x\}$ induces a $2K_2\cup K_1$, a contradiction. This proves \ref{C7-neiA}. $\Diamond$

 \begin{pro}\label{C7-x-Ai-2Ai+2}
Let $x\in V(H)\sm (A\cup D)$. Suppose $x$ has a neighbor in $A_i$. Then exactly one of  $N(x)\cap A_{i-2}$, $N(x)\cap A_{i+2}$ is nonempty.
\end{pro}
\no{\it Proof of $\ref{C7-x-Ai-2Ai+2}$}.~Suppose not, and let $i=1$. Let $a$ be a neighbor of $x$ in $A_1$. If $N(x)\cap A_{3}=\es$ and $N(x)\cap A_{6}=\es$, then by  \ref{c5-free}, $N(x)\cap A_{5}=\es$, and then  $\{a,x,u_5,u_6,u_3\}$ induces a $2K_2\cup K_1$, a contradiction; so we may assume that $N(x)\cap A_{3}\neq\es$  and $N(x)\cap A_{6}\neq\es$.  Then by \ref{c5-free}, $x$ is complete to $A_4\cup A_5$.  Then, again by using \ref{c5-free}, we see that $x$ is complete to $A_2\cup A_7$.  But then $x\in D$, a contradiction. This proves \ref{C7-x-Ai-2Ai+2}. $\Diamond$

 \begin{pro}\label{C7-Part}
 $V(H)=A\cup B \cup D$.
 \end{pro}
\no{\it Proof of $\ref{C7-Part}$}.~Let $x\in V(H)\sm (A\cup D)$.  Then, by \ref{C7-neiA}, we may assume that $x$ has a neighbor in $A_i$, say $a_{i}$. By \ref{C7-x-Ai-2Ai+2}, we may assume that   $N(x)\cap A_{i+2}\neq\es$ and $x$ is anticomplete to $A_{i-2}$. Then, by  \ref{c5-free}, $x$ is anticomplete  $A_{i-3}$. Let  $a_{i+2}$ be a neighbor of $x$ in $A_{i+2}$.  We claim that $x$ has a neighbor in $A_{i+1}$. Suppose  $x$ is anticomplete to $A_{i+1}$. Then, by \ref{c5-free}, $x$ is anticomplete to $A_{i+3}\cup A_{i-1}$. Also, if $x$ has a nonneighbor, say $a_i'$, in $A_{i}$, then $\{a_i',u_{i-1},x,a_{i+2},u_{i-3}\}$ induces a $2K_2\cup K_1$, a contradiction; so $x$ is complete to $A_i$. Likewise,  $x$ is complete to $A_{i+2}$. But then $x$ can be added
to $A_i$, contradicting the maximality of $A$. So we may assume that $x$ has a neighbor in $A_{i+1}$, say $a_{i+1}$. Then by \ref{c5-free}, $x$ has no neighbors in both  $A_{i+3}$ and $A_{i-1}$. But since $\{x,a_{i+1},u_{i+3},u_{i-3},u_{i-1}\}$ does not induce a $2K_2\cup K_1$, $x$ has a neighbor in  exactly one of $A_{i+3}$ and $A_{i-1}$, say $x$ has a neighbor in $A_{i+3}$.  So $x\in B_i$.
 This proves \ref{C7-Part}. $\Diamond$

 \begin{pro}\label{C7-Ai-stable}
 $A_i$ is a stable set.
 \end{pro}
 \no{\it Proof of $\ref{C7-Ai-stable}$}.~If there are adjacent vertices in $A_i$, say $p$ and $q$, then $\{p,q,u_{i+2},u_{i+3},u_{i-2}\}$ induces a $2K_2\cup K_1$, a contradiction. This proves \ref{C7-Ai-stable}. $\Diamond$

 \begin{pro}\label{C7-wi}
 $H[B_i]$ is $K_2\cup K_1$-free.
 \end{pro}
\no{\it Proof of $\ref{C7-wi}$}.~If there is a $K_2\cup K_1$ induced by the vertices, say $\{p,q,r\}$, in $B_i$, then $\{u_{i-1}, u_{i-2},p,q,r\}$ induces a $2K_2\cup K_1$, a contradiction. This proves \ref{C7-wi}. $\Diamond$

 \begin{pro}\label{C7-Bi-nei-comp}
$B_i$ is complete to $A_i\cup A_{i+1}\cup A_{i+2}\cup A_{i+3}$.
\end{pro}
\no{\it Proof of $\ref{C7-Bi-nei-comp}$}.Let $x\in B_i$ and $y\in A_i\cup A_{i+1}\cup A_{i+2}\cup A_{i+3}$, and suppose $x,y$ are nonadjacent. Let $a_{i+1}$ and $a_{i+2}$ be  neighbors of $x$ in $A_{i+1}$ and $A_{i+2}$ respectively. By symmetry, we may assume that $y\in A_i\cup A_{i+1}$.
 Now if $y\in A_i$, then $\{u_{i-1},y,x,a_{i+2},u_{i-3}\}$ induces a $2K_2\cup K_1$, a contradiction, and if $y\in A_{i+1}$, then, by \ref{C7-Ai-stable}, $ya_{i+1}\notin E$, and then $\{u_{i-1},u_{i-2},x,a_{i+1},y\}$ induces a $2K_2\cup K_1$, a contradiction.   This proves \ref{C7-Bi-nei-comp}. $\Diamond$

\begin{pro} \label{C7-BiBi+1}
$B_i$ is complete to $B_{i+1}\cup B_{i-1}$.
\end{pro}
\no{\it Proof of $\ref{C7-BiBi+1}$}.~Let $x\in B_i$ and $y\in B_{i+1}\cup B_{i-1}$, and suppose $x$, $y$ are nonadjacent. By symmetry, we may assume that $y\in B_{i+1}$. Then by \ref{C7-Bi-nei-comp}, $\{x, u_{i+1}, y, u_{i+4}, u_{i+3}\}$ induces a $P_5^c$, a contradiction. This proves \ref{C7-BiBi+1}. $\Diamond$

 \begin{pro}\label{C7-Bi-nemp}
If $B_i\neq\es$, then $B_{i-3}\cup B_{i-2}\cup B_{i+2}\cup B_{i+3}$ is empty.
 \end{pro}
\no{\it Proof of $\ref{C7-Bi-nemp}$}.~Let $x\in B_i$. Suppose  that there is a vertex $x'\in B_{i+2}$. If $xx'\in E$, then by \ref{C7-Bi-nei-comp}, $x$-$u_i$-$u_{i-1}$-$u_{i-2}$-$x'$-$x$ is a $C_5$, a contradiction; so $xx'\notin E$, and then, by \ref{C7-Bi-nei-comp}, $\{x,u_{i+1},x',u_{i-3},u_{i-1}\}$ induces a $2K_2\cup K_1$, a contradiction. So $B_{i+2}=\es$. Likewise, $B_{i-2}=\es$.
 Also, if there is a vertex, say $y\in B_{i+3}$, then, by \ref{C7-Bi-nei-comp},   $\{u_{i-1},u_{i},x,u_{i+3},y\}$   induces a $C_5$ or a $P_5^c$, a contradiction. So $B_{i+3}=\es$. Likewise, $B_{i-3}=\es$. This proves \ref{C7-Bi-nemp}. $\Diamond$

 \begin{pro} \label{C7-DUni}
 $D$ is complete to $A \cup B$.
 \end{pro}
 \no{\it Proof of $\ref{C7-DUni}$}.~Suppose there are nonadjacent vertices, say $x\in D$ and $a\in A_i$. Pick neighbors of $x$ in each  $A_{i+1}$, $A_{i+2}$ and $A_{i-1}$, say $p,q,$ and $r$ respectively. Then $\{a,p,q,r,x\}$ induces a $P_5^c$, a contradiction. So $D$ is complete to $A$. Next, if there are nonadjacent vertices, say  $x\in D$ and $x'\in B_i$, then, by \ref{C7-Bi-nei-comp}, and by the earlier argument, $\{x',u_i,u_{i-1},x,u_{i+3}\}$ induces a $P_5^c$, a contradiction.  This proves \ref{C7-DUni}. $\Diamond$

\smallskip

Now since $H^c$ is connected, we have $D=\es$. So by above properties, if $B=\es$, then $G$ is a clique-blowup of $C_7^c$. So we may assume that $B_1\neq\es$. Then by \ref{C7-Bi-nemp}, $B_3\cup B_4\cup B_5\cup B_6$ is empty, and one of $B_2$, $B_7$ is empty. Thus we conclude that $G$ is a blowup of $H^*$. Let $H^*$ be  defined as earlier. By the definition of blowup, $V(G)$ is partitioned into $Q_{v_i}$, $v_i \in V(H^*)$, such that each $Q_{v_i}$ induces a $P_3$-free graph. Now we let  $S_1:=R_{Q_{v_1}}\cup R_{Q_{v_4}}\cup R_{Q_{v_9}}$, $S_2:=R_{Q_{v_2}}\cup R_{Q_{v_5}}$, and $S_3:=R_{Q_{v_3}}\cup R_{Q_{v_7}}\cup R_{Q_{v_8}}$. Clearly $S_1,S_2$ and $S_3$ are  stable sets such that $S_1\cup S_2\cup S_3$ meets each maximal clique of $G$ twice. So $G$ is nice. This completes the proof of Theorem~\ref{thm-C7c}. \hfill{$\Box$}

\medskip
\no{\it Proof of Theorem~\ref{thm:struc-new}}.
Since each $k$-wheel, for $k\geq 6$ has an induced $P_5$, the proof of each of the  item in Theorem~\ref{thm:struc-new} follows from Theorems \ref{thm-5wheel}, \ref{thm:P5wheel-free-C5} and \ref{thm-C7c}   respectively. \hfill{$\Box$}

\bigskip
\no{\bf Acknowledgements}.~This paper is dedicated to the memory of Professor Fr\'ed\'eric Maffray on his third death anniversary. The authors thank an anonymous referee for his/her helpful suggestions and comments.

\end{document}